\documentclass{amsart}       
\usepackage[margin=1.2in]{geometry}
\usepackage{hyperref}
\usepackage{graphicx}
\usepackage{fixltx2e}                    
\usepackage[utf8]{inputenc}            
\usepackage{amssymb,latexsym, dsfont, amsfonts, amsbsy, amsmath, mathrsfs, dsfont}            
\usepackage{mathtools}                   
\usepackage{caption}
\usepackage{bm}                          
\usepackage{enumerate}                   
\usepackage{verbatim}                    
\usepackage{url}                         
\usepackage{tabularx}
\usepackage{microtype}  
\usepackage[all, knot]{xy}
        \xyoption{arc} 
        \xyoption{web}                 
\makeatletter                            
\def\MT@register@subst@font{\MT@exp@one@n\MT@in@clist\font@name\MT@font@list
 \ifMT@inlist@\else\xdef\MT@font@list{\MT@font@list\font@name,}\fi}
\makeatother

\usepackage{color}

\usepackage{scalerel,stackengine}
\stackMath
\newcommand\what[1]{%
\savestack{\tmpbox}{\stretchto{%
  \scaleto{%
    \scalerel*[\widthof{\ensuremath{#1}}]{\kern-.6pt\bigwedge\kern-.6pt}%
    {\rule[-\textheight/2]{1ex}{\textheight}}
  }{\textheight}%
}{0.5ex}}%
\stackon[1pt]{#1}{\tmpbox}%
}

\newcommand{\const}[1]{\overline{#1}}
\newcommand{\termDef}[1]{\textbf{\textit{#1}}}

\newcommand{\?}{\ensuremath{\mkern0.4\thinmuskip}}   

\let\leq=\leqslant
\let\geq=\geqslant
\let\Box=\square                            

\newcommand{\N}{\mathbb{N}}

\newcommand{\height}{\mathbf{\mathtt{h}}}

\newcommand{\sL}{\text{\scriptsize{\L}}}

\def\SF{\textit{SFm}}
\def\PSF{\textit{PSFm}}

\let\alg=\bm                                    
\let\aclass=\mathbb                             
\let\class=\mathbb                             
\let\mod=\mathfrak							
	
\let\alg=\mathbf
\let\log=\mathcal

\let\epsilon=\varepsilon
\let\Lambda\varLambda
\let\Gamma\varGamma
\let\Delta\varDelta
\let\Lambda\varLambda
\let\Omega\varOmega
\let\Theta\varTheta
\let\Xi\varXi
\let\Pi\varPi
\let\Sigma\varSigma

\newcommand{\aFL}{\mathbf{FL_{ew}}}                   
\newcommand{\conc}{\hspace{-0.12cm}\parallel\hspace{-0.1cm}}
\newcommand{\lu}{\scriptsize{\L}}

\renewcommand\conc {{_\smile}}

\theoremstyle{plain}
\newtheorem{theorem}{Theorem}[section]
\newtheorem{lemma}[theorem]{Lemma}
\newtheorem{corollary}[theorem]{Corollary}

\newtheorem*{obs}{Observation}
\theoremstyle{definition}

\newtheorem{definition}[theorem]{Definition}
\theoremstyle{remark}

\begin{document}
\title{Undecidability and non-axiomatizability of modal many-valued logics}

\author{Amanda Vidal}

\maketitle

\vspace{-0.8cm}
\begin{center}
Artificial Intelligence Research Institute (IIIA - CSIC) \\
Campus UAB, 08193 Bellaterra, Spain \\
\textsf{amanda@iiia.csic.es}
\end{center}
\vspace{0.2cm}

\thispagestyle{empty}

\begin{abstract}
	In this work we  study the decidability of a class of global modal logics arising from Kripke frames evaluated over certain residuated lattices, known in the literature as modal many-valued logics. We exhibit a large family of these modal logics which are undecidable, in contrast with classical modal logic and propositional logics defined over the same classes of algebras. This family includes the global modal logics arising from Kripke frames evaluated over the standard \L ukasiewicz and Product algebras. We later refine the previous result, and prove that global modal \L ukasiewicz and Product logics are not even recursively axiomatizable. 
	 We conclude by closing negatively the open question of whether each global modal logic coincides with its local modal logic closed under the unrestricted necessitation rule.
\end{abstract}

\section{Introduction}
Modal logics are one of the most developed and studied families of non-classical logics, exhibiting a beautiful equilibrium between complexity and expressivity. Generalizations of the concepts of necessity and possibility offer a rich setting to model and study notions from many different areas, including provability predicates, temporal and epistemic concepts, work-flows in software applications, etc. On the other hand, many-valued logics
provide a formal framework to manage gradual and resource sensitive information in a very general and adaptable way. 
In this work, we will focus our attention in the so-called  continuous t-norm fuzzy logics, a family of many-valued logics that has received special interest due to, among other properties, their completeness with respect to algebras on the real unit interval $[0,1]$ and their ability to capture the natural order with the implication operation.
Modal many-valued logics lie at the intersection of both modal and many-valued logics, extending many-valued logics with modal-like operators. 

This paper contributes to the problems of axiomatizability and decidability in these logics. We first show that global modal \L ukasiewicz (\L) and Product ($\Pi$) logics are undecidable, and we later refine those results to prove that they do not belong to $\Sigma_1$. Since the modal logics we investigate can be seen as fragments of the corresponding first order (F.O.) many-valued logics, our results likewise affect these fragments and the corresponding model theory \cite{CiH06}. 
On the other hand, 
the categorical equivalence between MV algebras (the algebraic semantics of \L ukasiewicz logic) and abelian lattice-ordered groups with a strong unit (\cite{Ch59, CiMuOt99}) directly relate the results shown here to the theory of $\ell$-groups.  Moreover, since the propositional basis of continuous logics \cite{ChKe66,YaPe10} (namely, the logical system proposed to study the so-called continuous model theory) coincides with \L ukasiewicz logic with an additional definable constant, our results also apply to modal continuous logics \cite{Ba18,Ba19} and to fragments of F.O. continuous logics \cite{YaPe10,YaUs10}, with the natural ramifications towards continuous model theory. Lastly, the field of modal many-valued logics is intrinsically related to that of fuzzy description logics (FDL) \cite{St01}, which also makes the results presented in this paper applicable to certain FDLs over the standard \L ukasiewicz and Product algebras (namely, over the real unit interval $[0,1]$).

 The notion of modal many-valued logic studied in this paper follows the tradition initiated by Fitting \cite{Fi92a,Fi92b} and Hájek \cite{HaHa96, Ha98}, which  differs from another relevant definition of so-called modal substructural logics studied for instance in \cite{O93, Re93b, Kam02}. The logics studied in this work are defined over valued Kripke models: Kripke frames $\langle W, R \rangle$ enriched with a world-wise evaluation\footnote{In a more general setting, also the accessibility relation can be evaluated over the same algebra. Nevertheless, in this paper, unless stated otherwise, whenever we talk about some modal many-valued logic we will be referring to the case with classical Kripke frames.} of the formulas $e \colon W \times Fm \rightarrow \alg{A}$
 into some algebra $\alg{A}$, and where $\Box, \Diamond$ are unary operators generalizing those of classical modal logic. 
 The operators $\Box$ and $\Diamond$ are not inter-definable in general \cite{RoVi20} in the minimal logics, and
 the logics with the two modalities might not be axiomatized by the addition of the axiomatic systems of the mono-modal (i.e., with only the $\Box$ or the $\Diamond$ operator) fragments. Their respective computational behavior might also differ. The previous facts imply that each one of these fragments needs to be, in general, separately studied.
 On the other hand, the so-called local and global logics' derivations arising from these models refer to the interpretation of the premises and conclusion in the derivability relation. In the former case they are considered world-wise, while in the latter the premises should hold in the whole model. These two semantics
 behave with respect to the F.O. semantics of the corresponding many-valued logic in the analogous way to the classical case and in Fisher-Servi Intuitionistic modal logic \cite{FS77}.

We will begin by addressing the question of the decidability of modal many-valued logics for minimal logics (namely, those that do not restrict the class of models). It is known that, in contrast to F.O. logic, which is undecidable, the minimal (classical) modal logic $K$ is, as propositional logic, decidable. 
In many-valued logics, similarly, F.O. logics are in most cases undecidable, while the propositional cases are usually co-NP complete. Nevertheless, while it is known that classical F.O. logic is $\Sigma_1$-complete, tautologies of F.O. over the standard \L ukasiewicz algebra form a $\Pi_2$ complete set \cite{Ra83, Ra83c}, and those over the standard Product algebra are $\Pi_2$-hard \cite{Ha98}. 
Regarding modal many-valued logics, the known results about decidability are rather partial.
G\"odel modal logics do not enjoy in general the finite model property with respect to the intended semantics \cite{CaRo10}. Interestingly enough,  in \cite{CaMe13, CaMeRo17, BaDi21}, the decidability of the local
 consequence relation for the minimal mono-modal and bi-modal logics, as well as for the S4 and S5 extensions, is proven. However, regarding the ongoing work, we mention that the decidability of the global consequence over the previous classes of models is still an open problem. 
 It is also known that the minimal local modal (standard) \L ukasiewicz and Product logics are decidable \cite{Vi20,Vi22}. On the other hand, in \cite{Vi20}, the undecidability of the local deductions over transitive models, valued respectively over the standard \L ukasiewicz and Product algebras is shown. 
Some additional decidability results have been proven in the context of modal many-valued logics where the accessibility relation of the Kripke models is also many-valued over the corresponding algebra. These results arise indirectly from the studies over Fuzzy Description Logics, which can be roughly interpreted as a fragment of multi-modal logic with many-valued accessibility and explicit rational constants, and where additional connectives are often included. Concerning the pure \L ukasiewicz and Product FDL, in \cite{BaPe11,BaPe11a, BoDiPe15, CeEs22} (see \cite{CeEs18} for a presentation of the previous results expressed in a modal logic setting) it is proven that the $r$-satisfiability question (for $r$ rational) is undecidable for the global modal Product and \L ukasiewicz cases, while the local $r$-satisfability problems are decidable. Nevertheless, the problems of validity and logical entailment remains open also for these modal logics with many-valued accessibility. More pertinently for the present paper, the fact that the accessibility relation is many-valued seems to be unavoidable in all the above proofs, which has shed little light on the setting based on classical Kripke frames.\footnote{For the interested reader, it is worth pointing out that the undecidability proofs from the previous references rely on infinite models, which in turn do not allow the use of those results to prove non-recursive enumerability results for those logics.}

In this paper, we study the decidability of entailment from finite sets of premises  in global modal logics. We show this question is undecidable for a large class of modal logics whose algebras of evaluation satisfy certain basic conditions (Theorem \ref{thm:Psatiff}). This class includes the minimal (bi)modal logics over the standard \L ukasiewicz and Product algebras.  The main problem that remains open concerning the decidability of minimal modal logics based on continuous t-norm logics is that of the global bi-modal G\"odel logic. 

Since the previous undecidability results are proven using classes of finite models, it will be possible to refine them to answer open problems concerning the axiomatizability of the same logics.
Let us briefly overview the known results and open questions concerning axiomatizability of the minimal modal fuzzy logics. 
In \cite{CaRo10, CaRo15, MeOl11, RoVi20} all the minimal modal logics associated to the standard G\"odel algebra are axiomatized. This includes both local and global deductions of the two mono-modal fragments and of the bi-modal logic with both $\Box$ and $\Diamond$.
In \cite{BoEsGoRo11}, a general study is done of the logics with only $\Box$ arising from 
Kripke models valued over finite residuated lattices, including axiomatization of the local and global deductions.
The axiomatic systems proposed there rely on the addition of canonical constants\footnote{Namely, one constant symbol for each element of the propositional algebra.} which among other things, make $\Box$ and $\Diamond$ interdefinable in modal logics arising from Kripke models valued over finite algebras \cite{Vi20a}.

On the other hand, regarding modal \L ukasiewicz (and the analogous continuous logic) and Product logics, the question of axiomatizability has not received a conclusive answer in the literature. The known results have proposed axiomatic systems that include some infinitary inference rule (i.e., with infinitely many premises), which are complete with respect to the infinitary deductions of the corresponding local or global logic. This solution is however non conclusive, since there is no explicit axiomatization for any of the intended finitary consequence relations, and questions related to the complexity  of the logics cannot be tackled as usual with an infinitary axiomatic system.\footnote{E.g., if a deductive system has a R.E. axiomatization with no infinitary rules, it is at least R.E. This reasoning cannot be done if the axiomatic system is infinitary in the above sense.}
Regarding the (standard) \L ukasiewicz case\footnote{By this we refer to the modal logics arising from models evaluated over the standard MV algebra.}, we find in \cite{HaTe13} an axiomatization of the infinitary local and global modal logics with an infinitary inference rule. Since the \L ukasiewicz negation is involutive, in this case the $\Box$ and $\Diamond$ modal operators are inter-definable, and all minimal logics (in the sense of modal operators) coincide.
%
For continuous modal logic an axiomatic system essentially equivalent to that in \cite{HaTe13} is proposed in \cite{Ba18}, still requiring an infinitary rule.
A similar situation happens in subsequent works concerning the study of the (standard) modal Product logic \cite{ViEsGo16} and other infinite linearly ordered residuated lattices \cite{Vid-PhD},
 where the proposed axiomatic systems require the extension of the logic with a dense countable set of constants, and an infinitary inference rule that quantifies over all the previous constants in the language. The cases studied in \cite{ViEsGo16, Vid-PhD} concern the logic with both $\Box$ and $\Diamond$ modalities, and to the best of our knowledge no works have studied the mono-modal fragments of these logics. Interestingly enough, in all the previous logics the infinitary rules are purely propositional, and result in the infinitary completeness of the corresponding propositional logic with respect to their corresponding algebraic semantics (\cite{Ha98, ViBoEsGo17}).

By proving that global modal \L ukasiewicz and Product logics are not recursively enumerable we contribute to answering the previous open questions, since this implies these logics are not axiomatizable by a R.E. (finitary) axiomatization.
These results can be seen in relation to the celebrated result by Scarpellini \cite{Scar62} that states that the set of tautologies of the infinitely-valued F.O. \L ukasiewicz logic is not recursively enumerable (later refined in \cite{Ra83, Ra83c}, proving that it is  $\Pi_2$-complete). Similarly, in \cite{Ra83}, Ragaz also proved that the satisfiability problem for the monadic fragment of F.O. \L ukasiewicz logic is $\Pi_1$ complete and undecidable in the presence of at least four symbols, and Bou showed\footnote{In an unpublished work presented in the LATD 2012 Tutorial, see \cite{Bou12}.} that validity in this fragment with at least two symbols is undecidable. The non-axiomatizability results for modal logics we present in this paper, when translated to F.O. \L ukasiewicz and Product logics, imply that their respective two variable fragments are not recursively enumerable in the presence of at least three unary predicates and a binary one\footnote{More specifically, with a binary $\{0,1\}$-valued predicate.}.

In the last section of this paper, we study the relation between the local and the global modal deductions, particularly motivated by the peculiarities intrinsic to the \L ukasiewicz case: while the local deduction (and so, the set of tautologies) is decidable ($\Delta_1$), the global deduction is not recursively enumerable ($\Sigma_1$). We will see that, as a consequence, the global deduction cannot be axiomatized by the local one extended with the usual necessitation rule ($\varphi \vdash \Box \varphi$). This contrasts with all other known cases in the literature and allows us to answer negatively this open question,  raised in \cite{BoEsGoRo11}.

The paper is organized as follows. We start in Section \ref{sec:Prelim} by introducing some necessary preliminaries. In Section \ref{sec:undecGlobal} we study the decidability of a large family of residuated lattice-based global modal logics, and prove they are undecidable by reducing the Post Correspondence Problem to them. This includes the standard \L ukasiewicz and Product cases. In Section \ref{sec:consequences} we obtain  negative results concerning the axiomatization (in the usual finitary way) of some of the above logics, namely that the (finitary) global modal standard \L ukasiewicz and Product logics are not recursively enumerable. We conclude the paper in Section \ref{sec:necRule} by showing that a global modal logic might fail to be axiomatized by an axiomatization of its corresponding logic plus the necessitation rule (which holds true for all logics from Section \ref{sec:undecGlobal}).

\section{Preliminaries}\label{sec:Prelim} 

In this work, a logic is identified with a consequence relation \cite{Fo16}, as opposed to only a set of formulas. 
While the second approach is more common in the literature of modal logics \cite{ChZa97}, we opt for the former definition because the differences between local and global modal logics are lost if only the tautologies of the logic are considered. Observe that the lack of the Deduction Theorem in the global logic makes the implication and the logical consequence not interchangeable.

Given a set of variables $\mathcal{V}$ and an algebraic language $\mathtt{L}$, the set $Fm^{\mathtt{L}}(\mathcal{V})$ is the set of formulas built from $\mathcal{V}$ using the symbols from $\mathtt{L}$. Unless stated otherwise, $\mathcal{V}$ is a fixed denumerable set, and it will be omitted in the notation of the set of formulas, and if the language is clear from the context we will omit it as well. 
A \termDef{rule} in $Fm$ is a pair $\langle \Gamma, \varphi \rangle \in \mathcal{P}(Fm) \times Fm$. We say a rule is \termDef{finitary} whenever $\Gamma$ is a finite set.
A \termDef{logic} $\mathcal{L}$ over $Fm$ is a consequence relation on $Fm$, that is, a set of rules such that:
\begin{enumerate}
\item $\mathcal{L}$ is \termDef{reflexive}, i.e., for every $\Gamma \subseteq Fm$ and every $\gamma \in \Gamma$, $\langle \Gamma, \gamma\rangle \in \mathcal{L}$, 
\item $\mathcal{L}$ satisfies \termDef{cut}, i.e., if $\langle \Gamma, \phi \rangle \in \mathcal{L}$ for all $\phi \in \Phi$, and $\langle\Phi, \varphi \rangle \in \mathcal{L}$ then $\langle \Gamma, \varphi \rangle \in \mathcal{L}$,
\item $\mathcal{L}$ is \termDef{substitution invariant}, i.e., for each substitution $\sigma$, if $\langle \Gamma, \varphi \rangle \in \mathcal{L}$ then $\langle \sigma[\Gamma], \sigma(\varphi) \rangle \in \mathcal{L}$.
\end{enumerate}

Whenever $\langle \Gamma, \varphi \rangle \in \mathcal{L}$ we will write $\Gamma \vdash_{\mathcal{L}} \varphi$.
Given a set of rules $R$, we will write $R^l$ to denote the minimal logic containing the rules in $R$. We say that a set of rules $R$ axiomatizes a logic $\mathcal{L}$ whenever $R^l = \mathcal{L}$. Observe that in this sense, every logic is axiomatized at least by itself.

With computational questions in mind, we will be focused on logics determined by finitary rules.
A logic $\mathcal{L}$ is \termDef{finitary} whenever $\Gamma \vdash_{\mathcal{L}} \varphi$ if and only if $\Gamma_0 \vdash_{\mathcal{L}} \varphi$ for some finite $\Gamma_0 \subseteq_\omega \Gamma$.\footnote{As usual, $\subseteq_\omega$ denotes the finite subsethood relation.}
For convenience, we will denote by $\mathcal{L}^{fin}$ the set of finitary consequences of $\mathcal{L}$, namely $\mathcal{L}^{fin} \coloneqq \{\langle \Gamma, \varphi\rangle \in \mathcal{L}\colon \Gamma \subseteq_{\omega} Fm\}$.

As it is most usual, in this work we do not consider infinite inputs for the computability questions. 
We say that a logic $\mathcal{L}$ is \termDef{decidable, recursive}, or \termDef{recursively enumerable}, respectively, if this is the case for the set $\mathcal{L}^{fin}$.

When $R$ is a set of finitary rules, the logic $R^l$ can be equivalently characterized through the usual notion of finite proof\footnote{There exists also a more general notion of proof managing infinitary rules, based on wellfounded trees, that we will not use here.} in $R$. Given a finite set of formulas $\Gamma \cup \{\varphi\}$, a \termDef{proof} or \termDef{derivation} of $\varphi$ from $\Gamma$
in $R$ is a finite list of formulas $\psi_1, \ldots, \psi_n$ such that $\psi_n = \varphi$ and for each $\psi_i$ in the list, either $\psi_i \in \Gamma$ or there is a rule $\Sigma \vdash \phi$ in $R$ and a substitution $\sigma$ such that $\sigma(\phi) = \psi_i$ and $\sigma[\Sigma]$ (possibly empty) is a subset of $\{\psi_1, \ldots, \psi_{i-1}\}$ (or empty if $i = 1$).
It is well known that $R^l = \{\langle \Gamma, \varphi \rangle \colon \Gamma \cup \{\varphi\}\subseteq Fm \text{ and there is a proof of }\varphi \text{ from }\Gamma \text{ in } R\}$.

We will say that a logic is \termDef{axiomatizable} whenever it can be axiomatized using a recursive set of finitary rules. It is clear that an axiomatizable  logic is finitary and R.E. On the other hand, a finitary R.E. logic with a definable idempotent $n$-ary operation for every $n$ is always axiomatizable. This fact is a natural generalization of Craig's Theorem, and it can be checked in a similar way. For the interested reader we provide the details in the Appendix.

In the rest of the paper we will work with logics having such an idempotent operation (which will be simply $\wedge$). Thereby we will resort without further notice to the following observation.

\begin{obs}
A finitary logic is recursively enumerable if and only if it is axiomatizable.
\end{obs}

Modal many-valued logics arise from Kripke structures evaluated over certain algebras, putting together relational and algebraic semantics in a way adapted to model different reasoning notions. In the next section, the general algebraic setting of these semantics will be the one of $FL_{ew}$-algebras, the corresponding algebraic semantics of the Full Lambek Calculus with exchange and weakening. This will offer a very general approach to the topic while relying on well-known algebraic structures. We will later focus on modal expansions of MV and product algebras, particular classes of $FL_{ew}$-algebras.

\begin{definition}
	An \termDef{$\aFL$-algebra} is a structure $\alg{A} = \langle A;  \wedge, \vee, \cdot, \rightarrow, \const{0}, \const{1}\rangle$ such that
	\begin{itemize}
		\item $\langle A; \wedge, \vee, \const{0},\const{1}\rangle$ is a bounded lattice;
		\item $\langle A; \cdot, \const{1}\rangle$   is a commutative monoid;
		\item $a \cdot b \leq c$  if and only if $a \leq b \rightarrow c$  for all $a,b,c\in A$. 
	\end{itemize}	
\end{definition}
We will usually write $ab$ instead of $a \cdot b$, and abbreviate $\overbrace{x \cdot x \cdots x}^{n}$ by $x^n$ for $n \geq 1$. Moreover, as usual, we will define $\neg a$ to stand for $a \rightarrow \const{0}$. A chain is a linearly ordered $\aFL$-algebra.

In the setting of the previous definition, we will denote by $\alg{Fm}'$ the algebra of formulas built over a countable set of variables $\mathcal{V}$ using the language corresponding to the above class of algebras (i.e., $\langle \wedge/2, \vee/2, \cdot/2, \rightarrow/2, \neg/1, \const{0}/0, \const{1}/0\rangle$). We will refer to the bottom and top elements of the algebra, $\const{0}$ and $\const{1}$, simply by $0$ and $1$. Moreover, we will again write $\varphi \psi$ instead of $\varphi \cdot \psi$ and $\varphi^n$ for the product of $\varphi$ with itself $n$ times (for $n \geq 1$), and we let, as usual
\[(\varphi \leftrightarrow \psi) \coloneqq (\varphi \rightarrow \psi) \cdot (\psi \rightarrow \varphi) \quad \text{and} \quad \neg \varphi \coloneqq \varphi \rightarrow \const{0}.\] 

For a set of formulas $\Gamma \cup \{\varphi\}$ and a class of $\aFL$-algebras $\class{A}$, we write $\Gamma \models_{\class{A}} \varphi$ if and only if, for each $\alg{A} \in \class{A}$ and each $h \in Hom(\alg{Fm}',\alg{A})$, if $h(\gamma) = 1$ for each $\gamma \in \Gamma$, then $h(\varphi) = 1$ too.
We will write $\models_\alg{A}$ instead of $\models_{\{\alg{A}\}}$.
As expected, $\models_{\class{A}}$ is a logic, and by convenience we will write $\models_{\class{A}}$ instead of $\vdash_{\models_{\class{A}}}$.

$\class{FL}_{ew}$, the class of $\aFL$-algebras, is a variety studied in depth, see for instance \cite{On10}, \cite{GaJiKoOn07}.

Let us introduce some examples of well-known subvarieties of $\class{FL}_{ew}$.
\termDef{Heyting Algebras}, the algebraic counterpart of Intuitionistic logic, are $\aFL$-algebras where $\wedge = \cdot$. 
The variety of \termDef{G\"odel algebras}, $\aclass{G}$, (corresponding to intermediate G\"odel-Dummett logic $\log{G}$) is that of semilinear Heyting algebras, i.e., 
those satisfying $(a \rightarrow b) \vee (b \rightarrow a)=1$ for all $a,b$ in the algebra. \termDef{BL algebras}, the algebraic counterpart of 
Hájek Basic Logic $\log{BL}$, are semilinear $\aFL$ algebras where $a \cdot (a \rightarrow b) = a \wedge b$ for every $a, b$ in the algebra. 
The variety of \termDef{MV algebras} $\aclass{MV}$, algebraic counterpart of \L ukasiewicz logic $\log{\L}$, is formed by the involutive BL algebras (i.e., satisfying 
$\neg \neg a \leq a$), and that of \termDef{Product algebras} $\aclass{P}$ (corresponding to Product Logic $\Pi$),  is formed by those
BL algebras satisfying $\neg \neg a \leq (b \cdot a \rightarrow c \cdot a) \rightarrow (b \rightarrow c)$ and $a \wedge \neg a \leq 0$. 

Particular algebras in the previous classes are the so-called \termDef{standard} ones, whose universe is the standard unit real interval $[0,1]$ and the order (affecting the lattice $\wedge, \vee$ operations) is the standard one. Let us introduce explicitly the operations, which we will denote, for convenience, with the subscripts $G$, $\L$ and $\Pi$ (for G\" odel, \L ukasiewicz and Product logic respectively).

\begin{itemize}
		\item $[0,1]_G$, the \textbf{standard G\"odel algebra}, defines  
		\[a \cdot_G b \coloneqq a \wedge b \qquad \text{and} \qquad a \rightarrow_G b \coloneqq \begin{cases} 1 &\hbox{ if } a \leq b\\ b &\hbox{ otherwise} \end{cases}\]

	\item $[0,1]_\sL$, the \textbf{standard MV algebra}, defines 	
	\[a \cdot_\sL b \coloneqq \max\{0, a + b -1\} \qquad \text{and} \qquad a \rightarrow_\sL b \coloneqq \min\{1, 1 - a + b\}\]
	\item $\bf{MV_n}$, the \textbf{$n$-valued MV algebra} is the subalgebra of $[0,1]_\sL$ with universe $\{0, \frac{1}{n-1}, \ldots, \frac{n-1}{n-1}\}$.
	\item $[0,1]_\Pi$, the \textbf{standard Product ($\Pi$) algebra}, defines 
	\[ a \cdot_\Pi b \coloneqq a \times b \qquad \text{and} \qquad  a \rightarrow_\Pi b \coloneqq \begin{cases} 1 &\hbox{ if } a \leq b\\ b/a &\hbox{ otherwise} \end{cases}\] with $\times$ being the usual product between real numbers;
	
\end{itemize}
 It is known that the standard G\"odel, MV and Product algebras generate their corresponding varieties. They do so also as quasi-varieties, which implies the completeness of the logics (understood as consequence relations) with respect to the logical matrices over the respective standard algebra. In the case of G\"odel, it is also the case that the variety is generated as a generalized quasi-variety, while this fails for MV and Product algebras. The first claim amounts to saying that for each set of formulas $\Gamma \cup\{ \varphi\}$, it holds that $\Gamma \vdash_{\log{G}} \varphi$ if and only if $\Gamma \models_{[0,1]_G} \varphi$ (and if and only if $\Gamma \models_{\aclass{G}} \varphi$). For a finite set of formulas $\Gamma \cup \{\varphi\}$ it holds that 
 $\Gamma \vdash_{\sL} \varphi$ if and only if $\Gamma \models_{[0,1]_{\sL}} \varphi$ (if and only if $\Gamma \models_{\aclass{MV}} \varphi$); and 
$\Gamma \vdash_{\Pi} \varphi$ if and only if $\Gamma \models_{[0,1]_{\Pi}} \varphi$ (if and only if $\Gamma \models_{\aclass{P}} \varphi$). The last conditions might fail for infinite $\Gamma$.

Le us introduce some  other families of $\aFL$-algebras that will be of use later on.
	\begin{definition} Let $\alg{A}$ be an $\aFL$-algebra.
		\begin{itemize}
			\item $\alg{A}$ is \termDef{$n$-contractive} whenever $a^{n+1} = a^n$ for all $a \in A$.
			\item $\alg{A}$ is \termDef{weakly-saturated} if for any two elements $a,b \in A$, if $a \leq b^n$ for all $n \in \mathds{N}$ then $ab = a$.
		\end{itemize}
	\end{definition}

Observe that if $\alg{A}$ is n-contractive, the element $a^n$ is idempotent (namely $a^n\cdot a^n = a^n$) for every $a \in A$. Simple examples of these algebras include Heyting algebras ($1$-contractive), or $MV_n$ algebras ($(n-1)$-contractive). On the other hand, the  standard MV-algebra and product algebra are not $n$-contractive for any $n$.
Regarding weakly saturation, observe that if the element $\inf\{b^n\colon n \in \mathds{N}\}$ exists in a weakly saturated algebra, then  it is an idempotent element. Examples of weakly saturated algebras are the standard MV-algebra, the standard product algebra, as well as the algebras belonging to the generalised quasi-varieties generated by  them.

The  algebra of modal formulas $\alg{Fm}$ is built in the same way as $\alg{Fm'}$, expanding the language of $\aFL$-algebras with two unary operators $\Box$ and $\Diamond$. While it is clear how to lift an evaluation from the set of propositional variables $\mathcal{V}$ into an $FL_{ew}$-algebra to $\alg{Fm'}$, the semantic definition of the modal operators depends on the relational structure in  the following way.

\begin{definition}\label{def:AKmodel}
 Let $\alg{A}$ be an $\aFL$-algebra. An \termDef{$\alg{A}$-Kripke model} is a structure $\mod{M} = \langle W, R, e\rangle$ such that
 \begin{itemize}
 \item $\langle W, R\rangle$ is a Kripke frame.  That is to say, $W$ is a non-empty set of so-called worlds and $R\subseteq W \times W$ is a binary relation over $W$, called an accessibility relation. We will often write $Rvw$ instead of $\langle v,w \rangle \in R$;
 \item $e$ is a map from $W \times \mathcal{V}$ to $A$.
 	\end{itemize}
The evaluation $e$ of an $\alg{A}$-Kripke model $\mod{M} = \langle W, R, e\rangle$ is uniquely extended to a map from $W \times Fm$ to $A$ by letting:
\begin{align*}
e(v, \const{c}) \coloneqq& c \text{ for }c \in \{0,1\} & e(v, \varphi \star \psi) \coloneqq& e(v, \varphi) \star e(v, \psi) \text{ for }\star \in \{\wedge, \vee, \cdot, \rightarrow\}\\
e(v, \Box \varphi) \coloneqq& \bigwedge\limits_{\langle v,w\rangle \in R} e(w, \varphi) & e(v, \Diamond \varphi) \coloneqq& \bigvee\limits_{\langle v,w\rangle \in R} e(w, \varphi)
\end{align*} 

A model $\mod{M} = \langle W, R, e\rangle$ is \termDef{safe} whenever the values of $e(v, \Box \varphi)$ and $e(v, \Diamond \varphi)$ are defined for every formula $\varphi$ at each world $v \in W$. 
The class of \termDef{$\aFL$-Kripke models} is the (set) union of all safe $\alg{A}$-Kripke models, for $\alg{A} \in \aFL$. 
\end{definition}
We call a model $\mod{M}$ \termDef{directed} whenever there is some world $u \in W$ in it such that, for each $v \in W$, there is some path from $u$ to $v$ \footnote{Finite sequence of worlds $\langle w_0, w_1 \ldots, w_n\rangle$ with $n \in \mathds{N}$, $w_i \in W$ for all $1 \leq i \leq n$ and such that $u = w_0, Rw_iw_{i+1}$ and $w_n = v$.} in $\mod{M}$.

Regarding notation, given a class of models $\class{C}$, we denote by $\omega\class{C}$ the class of finite models in $\class{C}$ (namely, the models $\omega\class{C} \coloneqq \{\langle W, R, e\rangle \in \class{C} \colon \vert W \vert < \omega\}$\footnote{ Observe that the algebra of evaluation is not necessarily finite.}).
 On the other hand, for a class of algebras $\aclass{A}$ (or a single algebra $\alg{A}$) we write $K{\aclass{A}}$  ($K{\alg{A}}$) to denote the  class of safe Kripke models over the algebras in the class (or over the single algebra specified). Finally, in order to lighten the reading, we will let $K\L$ and $K{\Pi}$ to denote respectively $K{[0,1]_{\tiny{\L}}}$ and $K{[0,1]_\Pi}$.

Towards the definition of modal logics over $\aFL$-algebras relying on the notion of $\aFL$-Kripke models, it is natural to use the notion of truth world-wise being $\{1\}$ (in  order to obtain, world-wise, the propositional $\aclass{FL}_{ew}$ logic). With this in mind, for each $\alg{A}$-Kripke model $\mod{M}$  and $v \in W$  we say that $\mod{M}$ \termDef{satisfies a formula} $\varphi$ \termDef{in $v$}, and write $\mod{M},v \models \varphi$ whenever  $e(v, \varphi) = 1$. Similarly, we simply say that $\mod{M}$ \termDef{satisfies a formula} $\varphi$, and write $\mod{M} \models \varphi$ whenever for all $v \in W$ $\mod{M},v \models \varphi$. The same definitions apply to sets of formulas.

As in the classical case,  the previous definition of satisfiability gives rise to two different logics: the local logic and the global one. In this work we will focus on the study of the global logic, but in Section \ref{sec:necRule} we will point out some results involving the local modal logic as well.

\begin{definition} Let $\Gamma \cup \{ \varphi\} \subseteq_{\omega} Fm$, and $\class{C}$ be a class of safe $\aFL$-Kripke models.
	\begin{itemize}
			\item \termDef{$\varphi$ globally follows from $\Gamma$ in  $\class{C}$}, and we write $\Gamma \vdash_{\class{C}} \varphi$, whenever  for every $\mod{M} \in \class{C}$, \[\mod{M} \models \Gamma \text{ implies }\mod{M} \models \varphi.\]
		
		\item 
		\termDef{$\varphi$ locally follows from $\Gamma$ in  $\class{C}$}, and we write $\Gamma \vdash^l_{\class{C}} \varphi$, whenever for every $\mod{M} \in \class{C}$ and every $v \in W$, \[\mod{M}, v \models \Gamma \text{ implies }\mod{M}, v \models \varphi;\]
	\end{itemize}
If $\class{C}$ is clear from the context, we will simply write $\vdash$ and $\vdash^l$ instead.
\end{definition} 
For arbitrary $\Gamma \cup \{ \varphi\} \subseteq Fm$, we let
$ \Gamma \vdash_{\class{C}} \varphi \text{ whenever there is } \Gamma_0 \subseteq_\omega \Gamma \text{ such that } \Gamma_0 \vdash_{\class{C}} \varphi$,
and the analogous for the local logic.

For a single Kripke model $\mod{M}$, we write $\Gamma \vdash_{\mod{M}} \varphi$ instead of  $\Gamma \vdash_{\{\mod{M}\}} \varphi$. In a similar way, for a model $\mod{M}$ and a world $u \in W$ we write
	$\Gamma \not \vdash_{\langle \mod{M}, u \rangle} \varphi$ to denote that $\mod{M} \models \Gamma$ and $\mod{M},u \not \models \varphi$ (namely, $\varphi$ does not follow globally from $\Gamma$ in $\mod{M}$, and world $u$ witnesses this fact). 
		 In a more general setting, fixing a Kripke frame $\mod{F}$ and an algebra $\alg{A}$, we write $\Gamma \vdash_{\mod{F}_{\alg{A}}} \varphi$ whenever $\Gamma \vdash_{\mod{M}} \varphi$ for every safe $\alg{A}$-Kripke model $\mod{M}$ with underlying Kripke frame $\mod{F}$. Analogously, for a class of frames $\class{F}$ and a class of algebras $\class{C}$, we write $\Gamma \vdash_{\class{F}_{\class{C}}} \varphi$ whenever $\Gamma \vdash_{\mod{F}_{\alg{A}}} \varphi$ for each $\mod{F} \in \class{F}$ and each $\alg{A} \in \class{C}$.

Tautologies (formulas following from $\emptyset$) of $ \vdash^l_{\class{C}}$ and $\vdash_{\class{C}}$ coincide, and $\vdash^l_{\class{C}}$ is strictly weaker than $\vdash_{\class{C}}$, a trivial separating case being the usual necessitation rule $\varphi \vdash \Box \varphi$ (valid in the global case and not in the local one). Observe that $ \vdash^l_{\class{C}}$ and $\vdash_{\class{C}}$ are, by their definition, determined by the safe directed models generated from the models in $\class{C}$.

Also, the unraveling and filtration\footnote{Identifying worlds $v,w$ such that $e(v, \varphi) = e(w, \varphi)$ for every formula $\varphi$.} techniques can be applied to an arbitrary directed model, obtaining a directed tree that, from the logical point of view, behaves in its root as the original model (i.e., satisfies exactly the same global and local derivations). Even if the resulting tree might be infinite, all worlds in the tree are, by construction, at a finite distance from the root.
Thus, $\vdash_{K\class{C}} \ =\  \vdash_{K\class{C}^T}$, for $K\class{C}^T$ being the class of \termDef{safe directed trees generated by models in $K\class{C}$}. 

Some useful notions concerning Kripke models are the following ones.
	
	\begin{definition}\label{def:height}
		Given a Kripke model $\mod{M}$ and $w \in W$, we let the \termDef{height of $w$}, $\height(w)$ be the element in $\mathds{N}\cup \{\infty\}$
		\footnote{Where $x < \infty$ for each  $x \in \mathds{N}$.}
		given by
		\[\height(w) \coloneqq
		sup\{k \in \N:\?  \exists w_0, \ldots, w_k \text{ with } w_0 = w \text{ and } Rw_i w_{i+1} \text{ for all } 0 \leq i \leq k \}.\]
	\end{definition}
	
	Observe that if there exists some cycle in the model, all worlds involved in that cycle (and every predecessor of each world in that cycle in the model) have infinite height.
	
	\begin{definition}
		Let $\varphi$ be a formula of $Fm$. We let the \termDef{subformulas of $\varphi$} be
		the set inductively defined by
		\begin{eqnarray*}
	\SF(p) & \coloneqq & \{p\}, \text{ for $p$ propositional variable or constant}\\
\SF(\triangledown \varphi) & \coloneqq & \SF(\varphi) \cup \{\triangledown \varphi\} \hbox{ for } \triangledown \in \{\neg, \Box, \Diamond\} \\
\SF(\varphi_1 \star \varphi_2) & \coloneqq & \SF(\varphi_1) \cup \SF(\varphi_2) \cup \{ \varphi_1 \star \varphi_2\} \hbox{ for } \star \in \{\wedge, \vee, \cdot, \rightarrow \}
		\end{eqnarray*}
	We let the \termDef{propositional subformulas of $\varphi$} be the set inductively defined by
	\begin{eqnarray*}
		\PSF(p) & \coloneqq & \{p\}, \text{ for $p$ propositional variable or constant}\\
		\PSF(\triangledown \varphi) & \coloneqq & \{\triangledown  \varphi\} \hbox{ for } \triangledown  \in \{\Box, \Diamond\}\\
		\PSF(\neg \varphi) & \coloneqq & \SF(\varphi) \cup \{\neg \varphi\} \\
		\PSF(\varphi_1 \star \varphi_2) & \coloneqq & \SF(\varphi_1) \cup \SF(\varphi_2) \cup \{ \varphi_1 \star \varphi_2\} \hbox{ for } \star \in \{\wedge, \vee, \cdot, \rightarrow \}
	\end{eqnarray*}
		For $\Gamma$ a set of formulas we let $\textit{(P)SFm}(\Gamma) \coloneqq \bigcup_{\gamma \in \Gamma} \textit{(P)SFm}(\gamma).$
	\end{definition}

Let us finish the preliminaries by stating a well-known undecidable problem, that will be used in the next sections to show undecidability of some of the modal logics introduced above. Recall that given two numbers  $\mathtt{ x}, \mathtt{  y}$ in base $s \in \mathds{N}$, their concatenation $\mathtt{ x_\smile   y}$ is given by $\mathtt{  x}s^{\parallel \mathtt{ y} \parallel} + \mathtt{ y}$, where $\parallel \mathtt{y}\parallel$ is the number of digits of $\mathtt{  y}$ in base $s$.
\begin{definition}[\textbf{Post Correspondence Problem (PCP)}]\label{def:PCP}
An instance $P$ of the PCP consists of a list $\langle \mathtt{x_1}, \mathtt{y_1} \rangle \dots \langle  \mathtt{x_n}, \mathtt{y_n}\rangle$ of pairs of numbers in some base $s \geq 2$.
 A solution for $P$ is a sequence of indexes $i_1, \dots, i_k$ with  $1 \leq i_j \leq n$ such that \[\mathtt{x_{i_1}}{_\smile} \dots _\smile  \mathtt{x_{i_k}}= \mathtt{y_{i_1}}{_\smile} \dots{_\smile} \mathtt{y_{i_k}}. \]
\end{definition}

The decision problem for PCP is, given a PCP instance, to decide whether such a solution exists or not. This question is undecidable \cite{Post46}.

\section{Undecidability of global modal logics}\label{sec:undecGlobal}
In this section, unless stated otherwise, we let $\aclass{A}$ be a class of weakly-saturated $\aFL$ chains
such that for every $n \in \mathds{N}$ there is some $\alg{A}_n \in \aclass{A}$ such that
$\alg{A}_n$ is non $n$-contractive. That is to say, there is some $a \in A_n$ such that 
$a^{n+1} < a^n.$

Examples of such classes of algebras are $\{[0,1]_{\lu}\}$ $\{MV_n\colon n \in \mathds{N}\}$ and $\{[0,1]_{\Pi}\}$. 
Natural examples of classes of algebras not satisfying the above conditions are $\{[0,1]_G\}$ and the variety generated by it, and the varieties of MV and product algebras (since these are not classes of weakly-saturated chains).

Let the class of frames $\class{L}$ be the isomorphic copies of the frames in the set 
\[\bigcup_{k \in \mathds{N}}\langle \{i \in \mathds{N}\colon i \leq k\}, \{\langle i, i+1\rangle\colon i \in \mathds{N}, i < k\} \rangle.\]
Namely, $\class{L}$ is given by the frames whose structure is isomorphic to the one depicted in Figure \ref{figureGlobal}, for any $k \in \mathds{N}$.
\begin{figure}
	\begin{displaymath}
		\xymatrix{ \underset{1}{\bullet} \ar[r] & \underset{2}{\bullet} \ar@{.}[rr] & &  \underset{k-1}{\bullet} \ar[r] & \underset{k}{\bullet}}
	\end{displaymath}
	\caption{Structure of the frames in $\class{L}$}\label{figureGlobal}
\end{figure}

\begin{theorem}\label{th:undec}
For any class of frames $\class{F}$ such that $\class{L} \subseteq \class{F}$ the logic $\vdash_{\class{F}_\class{A}}$ is undecidable. In particular, the logics $\vdash_{K{\aclass{A}}}$ and $\vdash_{\omega K{\aclass{A}}}$ are undecidable.

	More precisely, the three-variable fragments of the previous logics are undecidable.
\end{theorem}

The previous theorem follows as a direct consequence of the next result.
\begin{theorem}\label{thm:Psatiff}
	Let $P$ be an instance of the Post Correspondence Problem.
	 Then we can recursively define a set $\Gamma_P \cup \{\varphi_P\} \subseteq_{\omega} Fm$ in three variables from $P$ for which the following are equivalent:
	\begin{enumerate}
		\item $P$ is satisfiable;
		\item $\Gamma_P \not \vdash_{K{\aclass{A}}} \varphi_P$;
		\item $\Gamma_P \not \vdash_{\omega K{\aclass{A}}} \varphi_P$.
		\item $\Gamma_P \not \vdash_{\class{L}_{\class{A}}} \varphi_P$.
	\end{enumerate}
\end{theorem}
Trivially, $(4) \Rightarrow (3)$ and $(3) \Rightarrow (2)$ for any $\Gamma_P \cup \{\varphi_P\} \subseteq_{\omega} Fm$. In what remains of this section we will first show that $(1) \Rightarrow (4)$, and afterwards, that both $(2) \Rightarrow (4)$ and $(4) \Rightarrow (1)$. To this aim, let us begin by explicitly defining a suitable set of formulas $\Gamma_P \cup \{\varphi_P\}$.

For $P = \{\langle \mathtt{ x_1}, \mathtt{ y_1} \rangle \dots \langle  \mathtt{ x_n}, \mathtt{ y_n}\rangle\}$ list of pairs of numbers in base $s$, we let $\Gamma_P$ be
the set with formulas in variables $\mathcal{V} = \{x, y, z\}$:
\begin{enumerate}
	\item $\neg \Box \const{0} \rightarrow (\Box p \leftrightarrow \Diamond p)$ for each $p \in \mathcal{V}$; 
	 \item $\neg \Box \const{0} \rightarrow (z \leftrightarrow \Box z)$;
	 \item $\bigvee\limits_{1 \leq i \leq n} (x \leftrightarrow (\Box x)^{s^{\parallel\mathtt{ x_i}\parallel}} z^{\mathtt{ x_i}}) \land  (y \leftrightarrow (\Box y)^{s^{\parallel\mathtt{ y_i}\parallel}} z^{\mathtt{ y_i}})$.

\end{enumerate}

Finally, let $\varphi_P = (x \leftrightarrow y)^2 \rightarrow (x \rightarrow xz) \vee z.$

Roughly speaking, variables $x$ and $y$ will store information on the concatenation of the corresponding elements of the PCP, while $z$ will have a technical role. 

Given a solution of $P$, it is not hard to construct a finite model globally satisfying $\Gamma_P$ and not $\varphi_P$.
\begin{proof}\textit{(of Theorem \ref{thm:Psatiff},} $\mathit{(1) \Rightarrow (4)}$)\\
	Let $i_1, \dots, i_k$ be a solution for $P$, so $\mathtt{ x_{i_1}} \conc \dots \conc  \mathtt{ x_{i_k}} = \mathtt{ y_{i_1}} \conc \dots \conc \mathtt{ y_{i_k}} = r$ for some $r \in \mathds{N}$. Pick some non $r$-contractive algebra $\alg{A} \in \aclass{A}$ and $a \in A$ such that $a^{r+1} < a^r$, and define a finite $\alg{A}$-Kripke model $\mod{M}$ as follows:
	\begin{itemize}
		\item $W \coloneqq \{v_1, \dots v_k\}$;
		\item $R \coloneqq  \{\langle v_s, v_{s-1}\rangle \colon 2 \leq s \leq k\}$;
		\item For each $1 \leq j \leq k$ let \begin{itemize}
			\item $e(v_j, z) = a$ ;
			\item $e(v_j, x) = a^{{\mathtt{x_{i_1}}}\conc \dots \conc  {\mathtt{x_{i_j}} }}$;
			\item $e(v_j, y) = a^{\mathtt{ y_{i_1}} \conc \dots \conc  \mathtt{ y_{i_j}}}$;
		\end{itemize}	 
	\end{itemize}
The formula $\neg \Box 0$ is evaluated to $0$ in $v_1$, and to $1$ in all other worlds of the model. Thus, since $z$ is evaluated to the same value in all worlds of the model, and each world has exactly one successor except for $v_1$ (which has none), clearly the family of formulas in $(1)$ and in $(2)$ from $\Gamma_P$ are satisfied in all worlds of the model.

 To check that formula $(3)$ from $\Gamma_P$ is satisfied in all worlds of the model we reason by induction on the height of the world (Definition \ref{def:height}). For $v_1$ (with height equal to $0$), given that it does not have any successors, it is clear that 
\begin{align*}
e(v_1, (3)) =\ & \bigvee_{1\leq j \leq n} (e(v_1, x) \leftrightarrow e(v_1, z)^{\mathtt{x_j}}) \wedge  (e(v_1, y) \leftrightarrow e(v_1, z)^{\mathtt{y_j}})\\
=\ & \bigvee_{1\leq j \leq n} (a^{\mathtt{x_{i_1}}} \leftrightarrow a^{\mathtt{x_j}}) \wedge (a^{\mathtt{y_{i_1}}} \leftrightarrow a^{\mathtt{y_j}})\\
\geq\ & (a^{\mathtt{x_{i_1}}} \leftrightarrow a^{\mathtt{x_{i_1}}}) \wedge (a^{\mathtt{y_{i_1}}} \leftrightarrow a^{\mathtt{y_{i_1}}}) = 1
\end{align*} 
 
 For all other $v_r$ with $r >1$, recall that its only successor is $v_{r-1}$. Applying the definition of concatenation 
 , and the fact that for every $\alg{A} \in \aclass{A}$ and each $a \in A$ and $n,m\in \mathds{N}$, trivially $a^n a^m = a^{n+m}$ and $(a^n)^m = a^{nm}$, we can prove that 
 \begin{align*}
 e(v_r,(3)) =\ & \bigvee_{1\leq j \leq n} (e(v_r, x) \leftrightarrow e(v_r, \Box x)^{s^{\parallel\mathtt{ x_j}\parallel}} e(v_r, z)^{\mathtt{x_j}}) \wedge  (e(v_r, y) \leftrightarrow e(v_r, \Box y)^{s^{\parallel\mathtt{ y_j}\parallel}} e(v_r, z)^{\mathtt{y_j}}) \\
=\ & \bigvee_{1\leq j \leq n} (a^{{\mathtt{x_{i_1}}}\conc \dots \conc  {\mathtt{x_{i_r}} }} \leftrightarrow e(v_{r-1}, x)^{s^{\parallel\mathtt{ x_j}\parallel}}  a^{\mathtt{x_j}}) \wedge  (a^{{\mathtt{y_{i_1}}}\conc \dots \conc  {\mathtt{y_{i_r}} }} \leftrightarrow e(v_{r-1}, y)^{s^{\parallel\mathtt{ y_j}\parallel}}  a^{\mathtt{y_j}}) \\
 =\ & \bigvee_{1\leq j \leq n} (a^{{\mathtt{x_{i_1}}}\conc \dots \conc  {\mathtt{x_{i_r}} }} \leftrightarrow (a^{{\mathtt{x_{i_1}}}\conc \dots \conc  {\mathtt{x_{i_{r-1}}} }}) ^{s^{\parallel\mathtt{ x_j}\parallel}}  a^{\mathtt{x_j}}) \wedge  (a^{{\mathtt{y_{i_1}}}\conc \dots \conc  {\mathtt{y_{i_r}} }} \leftrightarrow (a^{{\mathtt{y_{i_1}}}\conc \dots \conc  {\mathtt{y_{i_{r-1}}} }})^{s^{\parallel\mathtt{ y_j}\parallel}}  a^{\mathtt{y_j}}) \\
 =\ & \bigvee_{1\leq j \leq n} (a^{{\mathtt{x_{i_1}}}\conc \dots \conc  {\mathtt{x_{i_r}} }} \leftrightarrow (a^{{\mathtt{x_{i_1}}}\conc \dots \conc  {\mathtt{x_{i_{r-1}}}}\conc {\mathtt{x_j}}}) \wedge (a^{{\mathtt{y_{i_1}}}\conc \dots \conc  {\mathtt{y_{i_r}} }} \leftrightarrow (a^{{\mathtt{y_{i_1}}}\conc \dots \conc  {\mathtt{y_{i_{r-1}}}}\conc {\mathtt{y_j}}})\\
 \geq\ &\ (a^{{\mathtt{x_{i_1}}}\conc \dots \conc  {\mathtt{x_{i_r}} }} \leftrightarrow (a^{{\mathtt{x_{i_1}}}\conc \dots \conc  {\mathtt{x_{i_{r-1}}}}\conc {\mathtt{x_{i_r}}}}) \wedge (a^{{\mathtt{y_{i_1}}}\conc \dots \conc  {\mathtt{y_{i_r}} }} \leftrightarrow (a^{{\mathtt{y_{i_1}}}\conc \dots \conc  {\mathtt{y_{i_{r-1}}}}\conc {\mathtt{y_{i_r}}}}) = 1
\end{align*} 

With the above, we have proven that $\mod{M} \models \Gamma_P$. 

On the other hand, since $i_1, \dots i_k$ was a solution for $P$, $e(v_k, x) = e(v_k, y)$. Moreover,  $e(v_k, z) = a < 1$, and $e(v_k, xz) = a^{r+1} < a^r = e(v_k, x)$, so $e(v_k, xz \rightarrow x) < 1$. This implies that $e(v_k, x \leftrightarrow y)^2 \rightarrow e(v_k, z) \vee e(v_k, xz \rightarrow x) < 1$, proving that $\Gamma_P \not \vdash_{\omega K{\aclass{A}}} \varphi_P$.	
\end{proof}

In order to prove the other implications of Theorem \ref{thm:Psatiff}, let us first show some technical characteristics of the models satisfying $\Gamma_P$ and not $\varphi_P$.

A first easy observation is that in every model satisfying $\Gamma_P$, the variable $z$ takes the same value in all connected worlds of the model. Relying on the completeness with respect to trees, we can prove that, in these models, $z$ is evaluated to the same value in the whole model.
\begin{lemma}\label{lemma:valuey}
	Let $\alg{A} \in \aclass{A}$, and $\mod{M} \in K\alg{A}^T$ with root $u$ be such that $\Gamma_P \not \vdash_{\langle \mod{M}, u\rangle}\varphi_P$. Then there is $\alpha_z \in A$ such that, for each world $v$ in the model,  $e(v,z) = \alpha_z$.
\end{lemma}
\begin{proof}
	Let $\alpha_z = e(u, z)$. It is easy to prove the lemma by induction on the distance of $v$ from $u$, which is always finite because $K\alg{A}^T$ is a class of directed trees.
	
	If $v = u$ then the claim follows trivially. Otherwise, assume that there are $w_0, w_1, \ldots, w_{k+1} \in W$ 
	with $w_0 = u, w_{k+1} = v$ and such that $R w_i w_{i+1}$ for all $0 \leq i \leq k$.
	Since $e(w_k, (1)) = e(w_k, (2))  = 1 $ and $R w_kw_{k+1}$, then we know
	\[
	e(w_k, \Box z) = e(w_k, \Diamond z)\qquad \text{ and }\qquad
	e(w_k, z) = e(w_k, \Box z)
	\]
	From the first equality we get that $e(v_1, z) = e(v_2,z)$ for all $v_1,v_2\in W$ such that $R w_k v_1$ and $R w_k v_2$. In particular, this yields that $e(w_k, \Box z) = e(w_{k+1}, z)$. Together with the second equality, it follows that $e(w_k, z) = e(w_{k+1}, z) = e(v,z)$. Applying the Induction Hypothesis, we conclude $e(u, z) = e(w_k, z) = e(v, z)$.
\end{proof}

The fact that algebras in $\aclass{A}$ are linearly ordered and weakly saturated allows us to also prove that such models can be assumed to have finite height.
\begin{lemma}\label{lemma:finiteHeightGlobal}
	Let $\alg{A} \in \aclass{A}$, and $\mod{M} \in K\alg{A}^T$ with root $u$ be such that $\Gamma_P \not \vdash_{\langle \mod{M}, u\rangle}\varphi_P$.
	Then $u$ has finite height. 
\end{lemma}
\begin{proof}
	From Lemma \ref{lemma:valuey} we know that  in each world $v$ of $\mod{M}$ it holds that $e(u, z) = e(v, z) = \alpha_z$. 
	Moreover, from $(3)$ in $\Gamma_P$ it follows that
	\[e(u,x) \leq \alpha_z^n \qquad \text{ for all }n \in \mathds{N} \text{ such that } n \leq \height(u)\]
	If $u$ was of infinite height, by weak saturation of $\alg{A}$, it would follow that $e(u,x) e(u,z) = e(u,x)$. However, since $e(u, \varphi_P) <1$, necessarily $e(u, xz) < e(u, x)$, and thus $u$ must be of finite height.
\end{proof}

As a corollary, we get that the values of $x$ and $y$ at each world are powers of $\alpha_z$. 
\begin{corollary}\label{corollary:alphay}
	Let $\alg{A} \in \aclass{A}$, and $\mod{M} \in K\alg{A}^T$ with root $u \in W$ be such that $\Gamma_P \not \vdash_{\langle \mod{M}, u\rangle}\varphi_P$. 
	Then for each $v \in W$ there are $a_v, b_v \in \mathds{N}$ such that  
	\[e(v, x) = \alpha_z^{a_v} \qquad \text{ and } \qquad e(v,y) = \alpha_z^{b_v}\]
	Moreover, if $\height(v) < \height(w)$ then $a_v < a_w$ and $b_v < b_w$.
	\end{corollary}
	\begin{proof}
	The first part easily follows by induction on the height of the model, from the previous lemma and formulas $(1)$ and $(3)$ in $\Gamma_P$. The second claim is immediate for the case when $Rvw$, since from $(3)$ implies that $e(v, x) \leq e(w, x) \alpha_z$ (and the same for variable $y$). For arbitrary $\height(v) < \height(w)$, this process is iterated. 
	\end{proof}
	
Another corollary can be proven after observing the way in which the implication behaves between powers of the same element in $FL_{ew}$ chains.
\begin{lemma}
Let $\alg{A} \in \aclass{A}$. For every $m > n \in \mathds{N}$ and every $a \in A$ such that $a^{m+1} < a^m$, it holds that $(a^n \rightarrow a^m)^2 \leq a$.
\end{lemma}
\begin{proof}
	If $n+1 < m$ (i.e., $m = n+1+k$ for some $k \geq 1$), we know that $a^{n+1} > a^m$: otherwise $a^m = a^{n+1+k}  = a^{n+1}$ implying that $a^{m+1} = a^{n+2} =  a^{n+1} = a^m$ too, which contradicts the assumptions. Thus,  $a^{n+1} \rightarrow a^m < 1$. By residuation, this is equivalent to $a \rightarrow (a^n \rightarrow a^m) < 1$, which implies $a >  a^n \rightarrow a^m$. In particular, the latter is also greater or equal than $( a^n \rightarrow a^m)^2$. \\
	Otherwise, necessarily $n+1 = m$. Let $a^n \rightarrow a^{n+1} = b$ for some $b \in A$. By residuation, $ba^n \leq a^{n+1}$, and so, $bba^n \leq ba^{n+1} \leq a^{n+2}$. Again by residuation, it follows that $b^2 \leq a^n \rightarrow a^{n+2}$. This is now an implication falling in the previous case (with $n+1 < m' = n+2$). Thus, we know that $a^n \rightarrow a^{n+2} < a$. We conclude that $b^2 \leq a$.
\end{proof}

	\begin{corollary}\label{corollary:eux=euy}
	Let $\alg{A} \in \aclass{A}$, and $\mod{M} \in K\alg{A}^T$ with root $u \in W$ be such that $\Gamma_P \not \vdash_{\langle \mod{M}, u\rangle}\varphi_P$.  Then $e(u, x) = e(u,y)$.
\end{corollary}
\begin{proof}
Corollary \ref{corollary:alphay} implies $e(u, x \leftrightarrow y) =\alpha_z^a \leftrightarrow \alpha_z^b$ for some $a,b \in \mathds{N}$. 
 From the previous lemma  we 
get that either $e(u, x \leftrightarrow y) = 1$ or $e(u, x \leftrightarrow y)^2 \leq \alpha_z$. Since the second condition implies $e(u, \varphi_P) = 1$, and this is false,  necessarily $e(u, x) = e(u, y)$. 
\end{proof}

%

We can now prove that if $\Gamma_P \not\vdash_{K{\aclass{A}}} \varphi_P$ then this happens in a model whose frame belongs to the class $\class{L}$, introduced at the beginning of this section (Figure \ref{figureGlobal}).


\begin{lemma}\label{lemma:completenessGlobal}
	\[\Gamma_P \vdash_{K{\aclass{A}}} \varphi_P \text{ if and only if } \Gamma_P \vdash_{\class{L}_{\aclass{A}}} \varphi_P\]
\end{lemma}
\begin{proof}
	Left-to-right direction is immediate since $\class{L}_{\aclass{A}} \subseteq K{\aclass{A}}$. 
	Concerning the right-to-left direction, assume $\Gamma_P \not \vdash_{K{\aclass{A}}} \varphi_P$. We know then there is a model $\mod{M} \in K{\aclass{A}}^T$ and $u \in W$ such that $\Gamma_P \not \vdash_{\langle \mod{M}, u\rangle}\varphi_P$.

We let the submodel $\what{\mod{M}}$ of $\mod{M}$  be defined with universe  $\{v_i\colon i \in \mathds{N}, v_i \in W\}$ such that 
	$v_1 \coloneqq u$ and for each $i \in \mathds{N}$, either $Rv_iv_{i+1}$ or $v_i$ has no successors in $\mod{M}$ and $v_{i+1} = v_i$. Namely, the universe of $\what{\mod{M}}$ is any chain of worlds from $\mod{M}$ beginning by $u$ and ending in a world with no successors.
	
	Define the model $\what{\mod{M}}$ by restricting to $\what{W}$ the accessibility relation and the evaluation from $\mod{M}$.
	From Lemma \ref{lemma:finiteHeightGlobal} we know $u$ has finite height in the original model, and so also $\what{\mod{M}}$ is finite (since by definition, for any $j > \height(w)$, $v_j = v_{\height(w)}$). Henceforth, by construction, its underlying frame is isomorphic to the frame from $\class{L}$ with universe $\{i \in \mathds{N} \colon i \leq \height(w)\}$. Namely, $\what{\mod{M}} \in \class{L}_{\aclass{A}}$.

	It remains to prove that, for every $\psi \in SFm(\Gamma_P\cup\{\varphi_P\})$ and each $v \in \what{W}$, it holds that $\what{e}(v, \psi) = e(v, \psi)$. 
	It is clear that restricting to a submodel does not change the value of propositional variables at each world, 
	i.e., for every $p \in \mathcal{V}$ (and thus, also for every propositional formula) and any $t \in \what{W}$
	it holds that $\what{e}(v,p) = e(v,p)$. For other formulas, we prove the analogous claim by induction on the formula and on the height of $v$ in $\what{\mod{M}}$.
	
The case for $\height(v)= 0 $ (i.e., there are no successors) is trivial to check, 
	since by construction, $v$ 
	does not have successors in $\mod{M}$ either. Thus, all formulas beginning with a modality contained in 
	$\SF(\Gamma_P \cup \{\varphi_P\})$ are evaluated (both in $\mod{M}$ and in $\what{\mod{M}}$) to either 
	$\const{1}$ ($\Box$) or $\const{0}$ ($\Diamond$). 
	Since the values of the propositional variables are not modified 
	by taking submodels, this concludes the proof of the step.
	
	Regarding the case for $\height(v) = n+1$ in $\what{\mod{M}}$, observe that 
	$v$ has successors both in $\mod{\mod{M}}$ and $\what{\mod{M}}$, so
	$\what{e}(v, \Box \const{0}) = e(v, \Box \const{0}) = 0$. 
	On the other hand, 
	$e(v, \Box p) = e(v, \Diamond p)$ for all $p \in \mathcal{V}$ 
	(from formulas in $(1)$), and so, in all successors of $v$ in $\mod{M}$, each variable $p$ 
	takes the same value, say $\alpha_p$. 
	Then, in particular, in the world $w$ chosen as the only successor of $v$ in the construction of $\what{\mod{M}}$, 
	it also holds that $e(w, p) = \what{e}(w, p) = \alpha_p$. 
	Since by construction of $\what{\mod{M}}$ the world $v$ has as only successor $w$, it holds that $\what{e}(v, \Box p) = \what{e}(v, \Diamond p) = \what{e}(w, p)$. Then, 
	$\what{e}(v, \Box p) = \what{e}(v, \Diamond p) = e(v, \Diamond p) = e(v, \Box p) = e(w, p) = \alpha_p$.

	The only formulas beginning with a modality appearing in $\SF(\Gamma_P \cup \{\varphi_P\})$
	are of the form $\Box \const{0}$, $\Box p$ and $\Diamond p$ for $p \in \mathcal{V}$. 
	Since the evaluation of all these formulas and of the propositional variables from $\mathcal{V}$ 
	in the world $v$ coincides in $\mod{M}$ and $\what{\mod{M}}$
	we conclude that the evaluation in $v$ of formulas built from these ones using propositional connectives is also preserved from 
	$\mod{M}$ to $\what{\mod{M}}$.
\end{proof}

At this point, it is possible to obtain a useful characterization of $x$ and $y$ in terms of 
$\alpha_z$ at each world of a model with a frame isomorphic to one in $\class{L}$ that satisfies $\Gamma_P$ and not $\varphi_P$ in its root. For convenience, in the next result we will invert the labeling of the worlds in the frames isomorphic to those in $\class{L}$, namely resorting to frames with structure $\langle \{u_i \colon i \in \mathds{N}, i \leq k\}, \{\langle u_{i+1},u_i\rangle\colon i \in \mathds{N}, i < k\}\rangle$ (trivially isomorphic to $\langle \{i \colon i \in \mathds{N}, i \leq k\}, \{\langle i, i+1\rangle\colon i \in \mathds{N}, i < k\}\rangle \in \class{L}$).
%
%

\begin{lemma}\label{lemma:charvw}
	Let $\mod{M} = \langle \{u_i \colon i \in \mathds{N}, i \leq k\}, \{\langle u_{i+1},u_i\rangle\colon i \in \mathds{N}, i < k\}, e \rangle$ be an $\aclass{A}$-Kripke model
	such that $\Gamma_P \not \vdash_{\langle \mod{M}, u_k\rangle}\varphi_P$.
	Then there exist $i_1, \ldots, i_k$ with $0 \leq  i_j \leq n$ for each $1 \leq j \leq k$,  such that for each $1 \leq j \leq k$,
	\[e(u_j, x) = \alpha_z^{\mathtt{ x_{i_1}}\conc \dots\conc  \mathtt{ x_{i_j}}} \qquad \text{ and } \qquad e(u_j, y) = \alpha_z^{\mathtt{ y_{i_1}}\conc \dots\conc \mathtt{ y_{i_j}}}.\]
	
	Moreover, for each $1 \leq j \leq k$,  \[e(u_j,x) = e(u_j,y) \text{ if and only if } \mathtt{ x_{i_1}}\conc \dots\conc  \mathtt{ x_{i_j}} = \mathtt{ y_{i_1}}\conc \dots\conc \mathtt{ y_{i_j}} .\]
	
\end{lemma}
\begin{proof}
	We will prove the first claim of the lemma by induction on $j$. The details are only given for the $x$ case, the other one is proven in the same way.
	
	For $j = 1$, $u_1$ does not have successors. From formula $(3)$ in $\Gamma_P$ (relying on the fact that the algebras in $\aclass{A}$ are chains) and Lemma \ref{lemma:valuey} it follows that there is $i_1 \in \{1, \ldots, n\}$ for which
	\[e(u_1, x) = e(u_1, \Box x)^{s^{\parallel{\mathtt{ x_{i_1}}}\parallel}} e(u_1, z)^{\mathtt{ x_{i_1}}} = 1^{s^{\parallel\mathtt{ x_{i_1}}\parallel}}\alpha_z^{\mathtt{ x_{i_1}}} = \alpha_z^{\mathtt{ x_{i_1}}}.\]
	
	For $j = r+1$, observe the only successor of $u_j$ in $\mod{M}$ is $u_r$. Then, from $(3)$ and Lemma \ref{lemma:valuey} it follows that there is $i_j \in \{1, \ldots, n\}$ for which
	\[e(u_j, x) = e(u_j, \Box x)^{s^{\parallel\mathtt{ x_{i_j}}\parallel}}  e(u_j, z)^{\mathtt{ x_{i_1}}}= e(u_r, x)^{s^{\parallel\mathtt{ x_{i_j}}\parallel}} \alpha_z^{\mathtt{ x_{i_j}}}\]
	By the Induction Hypothesis, the above value is equal to  $
	 (\alpha_z^{{\mathtt{ x_{i_1}}\conc \dots\conc  \mathtt{ x_{i_r}}}})^{s^{\parallel\mathtt{ x_{i_j}}\parallel}} \alpha_z^{\mathtt{ x_{i_j}}}$, and through simple properties of the monoidal operation, to $ (\alpha_z^{{\mathtt{ x_{i_1}}\conc \dots\conc  \mathtt{ x_{i_r}}}})^{s^{\parallel\mathtt{ x_{i_j}}\parallel} + \mathtt{ x_{i_j}}}$ as well.
	This value, by definition of the concatenation of numbers in base $s$, is exactly $\alpha_z^{{\mathtt{ x_{i_1}}\conc \dots\conc  \mathtt{ x_{i_j}}}}$, concluding the proof of the first claim.
	
	Concerning the second claim, assume towards a contradiction that there is $1 \leq j \leq k$ such that $\mathtt{ x_{i_1}}\conc \dots\conc  \mathtt{ x_{i_j}} \neq \mathtt{ y_{i_1}}\conc \dots\conc  \mathtt{ y_{i_j}}$ 
	 and 
	$e(u_j, x) = \alpha_z^{\mathtt{ x_{i_1}}\conc \dots\conc  \mathtt{ x_{i_j}}} =  \alpha_z^{\mathtt{ y_{i_1}}\conc \dots\conc  \mathtt{ y_{i_j}}} =  e(u_j,y)$.
	If $\mathtt{ x_{i_1}}\conc \dots\conc  \mathtt{ x_{i_j}} < \mathtt{ y_{i_1}}\conc \dots\conc  \mathtt{ y_{i_j}}$, 
	it follows that
	$\alpha_z^{\mathtt{ x_{i_1}}\conc \dots\conc  \mathtt{ x_{i_j}}} \alpha_z^n = \alpha_z^{\mathtt{ x_{i_1}}\conc \dots\conc  \mathtt{ x_{i_j}}} $ for every $n \geq 0$. Thus, in particular, from Corollary \ref{corollary:alphay}
	$e(u_k, x) = e(u_j, x)$, and also $e(u_k, x) \alpha_z = e(u_k, x) = e(u_j, x)$.
	However, $\mod{M}, u_k \not \models \varphi_P$ implies that $e(u_k, x)  \alpha_z < e(u_k,x)$, reaching a contradiction.

	The proof is analogous if
	$\mathtt{ x_{i_1}}\conc \dots\conc  \mathtt{ x_{i_j}} > \mathtt{ y_{i_1}}\conc \dots\conc  \mathtt{ y_{i_j}}$.\qedhere
	
\end{proof}

All the previous technical lemmas lead to a simple proof of Theorem \ref{thm:Psatiff}.
\begin{proof}\textit{(of Theorem \ref{thm:Psatiff},} $\mathit{(2) \Rightarrow (3) \Rightarrow (1)}$)\\
	Assume condition $\mathit{(2)}$ of the lemma, i.e. $\Gamma_P \not \vdash_{K{\aclass{A}}} \varphi_P$. Lemma \ref{lemma:completenessGlobal} implies there is a model $\mod{M} \in \what{K\aclass{A}}$ and $u \in W$ such that $\Gamma_P \not \vdash_{\langle \mod{M}, u_k\rangle}\varphi_P$. Since all models in $\what{K\aclass{A}}$ are finite, this proves point $(3)$. 
From here, from Corollary \ref{corollary:eux=euy} we know that $e(u, x) = e(u, y)$. Then, by Lemma \ref{lemma:charvw}, it follows 	that there exist indexes $i_1, \ldots, i_k$ in $\{1, \ldots, n\}$ for which 
	$\mathtt{ x_{i_1}}\conc \dots\conc \mathtt{ x_{i_k}}=  \mathtt{ y_{i_1}}\conc \dots\conc \mathtt{ y_{i_k}}$. This is a solution for the Post Correspondence Instance $(P)$, concluding the proof of $\mathit{(3) \Rightarrow (1)}$.\qedhere
	\end{proof}

Let us conclude this section by pointing out that the undecidability results in \cite{Vi20}, affecting the local deduction over transitive models evaluated over classes of algebras like $\aclass{A}$, are immediately translatable to the global logic. Consider the class of frames $\class{L}^4$ given by isomorphic copies of the frames in the set:
 \[\bigcup_{i \in \mathds{N}}\langle\{j \in \mathds{N}\colon j \leq i\}, <\rangle,\]
 where $<$ is the usual strict order relation in $\mathds{N}$ restricted in each case to the corresponding universe.
 \begin{corollary}[Of Theorem 3.1 from \cite{Vi20}]
 For any class of transitive frames $\class{T}$ containing an isomorphic copy of $\class{L}^4$, the logic $\vdash_{\class{T}_\class{A}}$ is undecidable. In particular, the logic $\vdash_{K4\aclass{A}}$ is undecidable, for $K4\aclass{A}$ being all safe transitive Kripke models evaluated over the algebras in $\aclass{A}$. 
 
  More precisely, the three-variable fragments of the previous logics are undecidable.
 \end{corollary}
 \begin{proof}
 It is straightforward, since for any finite $\Gamma \cup \{\varphi\} \subseteq Fm$ and any class of transitive models $\class{C}$, it is routine to check that
\[\Gamma \vdash_{\class{C}} \varphi \text{ if and only if } \Gamma, \Box \Gamma \vdash^l_{\class{C}} \varphi. \qedhere\]
 \end{proof}

\section{Non axiomatizability of modal \L ukasiewicz and Product logics}\label{sec:consequences}
The undecidability of the previous family of modal logics over finite models raises the question of their axiomatizaiblity. In particular, it was  an open problem how to axiomatize the finitary standard modal \L ukasiewicz logic (\cite{HaTe13},\cite{DiaMet18}) and standard modal Product logic (\cite{ViEsGo16}). None of the minimal modal logics over the previous standard algebras has been axiomatized in the literature, including the logics arising from crisp-accessibility frames (nor those from many-valued frames) with either one or both modal operators. In the previous references, infinitary axiomatic systems complete with respect to related (but different) deductive systems have been proposed. For instance, their corresponding infinitary companions (in some cases, over extended languages). 

We close this open problem for the standard \L ukasiewicz and Product logics with a negative answer: these logics are in fact not axiomatizable, since their respective sets of valid consequences are not recursively enumerable. We will devote this section to prove the previous claims. For that, 
three properties turn out to be crucial: undecidability of the global consequence over finite models of the class, decidability of the propositional logic and completeness of the global consequence  with respect to certain well-behaved models (in these cases, in terms of witnessing conditions). We will prove this negative result for the modal expansion of the standard \L ukasiewicz logic. Then, the analogous result will follow for the Product logic, relying on the known isomorphism between the standard MV-algebra and a certain Product algebra.

The first one of the above properties was proven in Section \ref{sec:undecGlobal}. Let us show how decidability of the underlying propositional logic $\models_{\aclass{A}}$ implies that the set $\{\langle \Gamma, \varphi \rangle\colon \Gamma \subseteq_\omega Fm, \varphi \in Fm, \Gamma \not \vdash_{\omega K{\aclass{A}}} \varphi\}$ is R.E., which will allow us to conclude there is no possible axiomatization for the logics of finite models over those classes of algebras.

We first see that the global consequence relation over a finite frame is decidable as long as the underlying propositional consequence relation is decidable too.

\begin{lemma}\label{lemma:finitelyManyModels}
	Let $\mod{F}$ be a finite frame, and $\aclass{A}$ a class of residuated lattices for which $ \models_{\aclass{A}}$ is decidable. 
	Then $\vdash_{\mod{F}_{\aclass{A}}}$ is decidable.
\end{lemma}
\begin{proof}

Let $\mod{F} = \langle W, R \rangle$ a Kripke frame, $v \in W$ a world in the frame, $\psi$ a modal formula with variables in a finite set $\mathcal{V}$, and $x$ not in $\mathcal{V}$. Consider the extended set of propositional variables
	\[\mathcal{V}^* \coloneqq \{p^v \colon p \in \mathcal{V}, v \in W\} \cup \{x^v_{\triangledown\varphi} \colon \triangledown \in \{\Box, \Diamond\}, \triangledown \varphi \in \SF(\psi), v \in W\}.\]
	
	We recursively define the non-modal formula $\langle \psi, v \rangle^*$ over $\mathcal{V}^*$ as follows:
	\begin{align*}
	    \langle c, v\rangle^* \coloneqq& c \text{ for } c \in \{\const{0},\const{1}\} & 
		\langle p, v\rangle^* \coloneqq& p^v \text{ for } p \in \mathcal{V}\\
		\langle \varphi \star \chi, v\rangle^* \coloneqq& \langle\varphi, v \rangle^* \star  \langle\chi, v \rangle^* \text{ for }\star \in \{\cdot, \rightarrow\} & 
		\langle \triangledown \varphi, v\rangle^* \coloneqq& x_{\triangledown\varphi}^v \text{ for }\triangledown \in \{\Box, \Diamond\}
	\end{align*}
	Then, for $\Sigma$ a set of formulas we let $\langle \Sigma, v\rangle^* \coloneqq \{\langle \sigma, v \rangle^*\colon \sigma \in \Sigma\}$, where set of original variables is $\mathcal{V} \coloneqq  \bigcup \{\mathcal{V}ars(\sigma) \colon \sigma \in \Sigma\}$.
	Moreover, consider the formulas\footnote{These are proper formulas because $W$ is a finite set.}
	 \[\delta_\Box^v(\psi) \coloneqq x_{\Box \psi}^v \leftrightarrow \bigwedge\limits_{w \in W:Rvw}\langle \psi, w\rangle^*\qquad \text{ and } \qquad
	\delta_\Diamond^v(\psi) \coloneqq x_{\Diamond \psi}^v \leftrightarrow \bigvee\limits_{w \in W:Rvw}\langle \psi, w\rangle^*.\]
	
	From those, define the set of formulas \[\Delta^v(\Gamma, \varphi) \coloneqq \{\delta_\Box^v(\psi) \colon \Box \psi \in \SF(\Gamma, \varphi)\} \cup \{\delta_\Diamond^v(\psi) \colon \Diamond \psi \in \SF(\Gamma, \varphi)\}.\]

	We will now prove that 
	\[\Gamma \vdash_{\mod{F}_{\aclass{A}}} \varphi \text{ if and only if } \{\langle \Gamma, v\rangle^*,  \Delta^v(\Gamma, \varphi) \colon v \in W\}
\models_{\aclass{A}} \bigwedge\limits_{v \in W} \langle \varphi, v\rangle^*,\]
		which implies the lemma.
	
	To prove the right-to-left direction, assume $\Gamma \not \vdash_{\mod{F}_{\aclass{A}}} \varphi$. Then there is $\alg{A} \in \aclass{A}$, and an $\alg{A}$-Kripke model over $\mod{F}$ in which $e(v, \Gamma) \subseteq \{1\}$ for all $v$ and $e(v_0, \varphi) < 1$ for some $v_0 \in W$. Consider then the mapping $h \colon \mathcal{V}^* \rightarrow A$ defined by $h(p^v) = e(v, p)$, $h(x_{\Box \psi}^v) = e(v, \Box \psi)$ and $h(x_{\Diamond \psi}^v) = e(v, \Diamond \psi)$. It is easy to see that the extension of this mapping to a homomorphism into $\alg{A}$ satisfies $h(\langle\psi, v\rangle^*) = e(v, \psi)$ for every $\psi \in SFm(\Gamma, \varphi)$. Thus, it satisfies the premises in the derivation at the right side, since $e(v, \Gamma) \subseteq \{1\}$ for all $v$ and by the semantical definition of $\Box$ and $\Diamond$ in the model. On the other hand, it does not satisfy the consequence, since $e(v_0, \varphi) < 1$.
	
	Conversely, given a propositional homomorphism $h$ over some algebra $\alg{A} \in \aclass{A}$ satisfying \[\{\langle \Gamma, v\rangle^*,  \Delta^v(\Gamma, \varphi) \colon v \in W\}\] and not satisfying $\bigwedge\limits_{v \in W} \langle \varphi, v\rangle^*$, we can consider the $\alg{A}$-Kripke model over $\mod{F}$ that lets $e(v, p) = h(p^v)$. Since $h(\Delta^v(\Gamma, \varphi)) = \{1\}$, then $e(v,\psi) = h(\langle\psi, v\rangle^*)$ for every $\psi \in SFm(\Gamma, \varphi)$, concluding the proof.
\end{proof}

\begin{corollary}\label{cor:decfixcard}
	Let $\aclass{A}$ be a class of residuated lattices for which $\models_{\aclass{A}}$ is decidable, and $j \in \mathds{N}$. 
	Then the problem of determining whether a formula $\varphi$ follows globally from a finite set of formulas $\Gamma$ in all $\aclass{A}$-models of cardinality $j$ (denoted by $\vdash_{jK{\aclass{A}}}$) is decidable.
\end{corollary}
\begin{proof}
	There is a finite number of frames of cardinality $j$, and so, for each one, we can run the decision procedure from the above lemma.
\end{proof}

Exhibiting a recursive procedure enumerating the elements not belonging to $\vdash_{\omega K{\aclass{A}}}$ is now easy.
\begin{lemma}\label{lemma:recomplement}
Let $\aclass{A}$ be a class of residuated lattices for which $\models_{\aclass{A}}$ is decidable. Then the set $\{\langle \Gamma, \varphi \rangle \in \mathcal{P}_\omega(Fm) \times Fm \colon \Gamma \not \vdash_{\omega K{\aclass{A}}} \varphi\}$ is recursively enumerable.
\end{lemma}
\begin{proof}

	Let us enumerate all pairs $\langle \Gamma, \varphi \rangle \in \mathcal{P}_\omega(Fm) \times Fm$, and initialize $P$ as the empty set. 
	Now, for each $i \in \mathds{N}$, store $\langle \Gamma_i, \varphi_i\rangle$ in $P$. Then check, for each $\langle \Gamma, \varphi \rangle \in P$ and for each $j \leq i$, whether $\Gamma\vdash_{jK{\aclass{A}}} \varphi$.  This is a finite amount (since $P$ is always finite) of decidable operations (from Corollary \ref{cor:decfixcard}), thus a decidable operation. Whenever the answer is negative, return that pair and continue.
	
	To see that the previous procedure enumerates exactly $\{\langle \Gamma, \varphi \rangle \in \mathcal{P}_\omega(Fm) \times Fm \colon \Gamma \not \vdash_{\omega K{\aclass{A}}} \varphi\}$, pick an arbitrary $\langle \Gamma_i, \varphi_i \rangle \in \mathcal{P}_\omega(Fm) \times Fm$ (according to the initial enumeration).

	First suppose that $\Gamma \not \vdash_{\omega K{\aclass{A}}} \varphi$, and that this happens in some $\alg{A}$-model of cardinality $j$ for some $\alg{A}\in \aclass{A}$. Then, at step $\max\{i, j\}$	
	it is checked whether $\Gamma_i \vdash_{jK{\aclass{A}}} \varphi_i$. By assumption we know this condition does not hold, so the answer is negative and the pair $\langle \Gamma_i, \varphi_i\rangle$ is returned.
	
	On the other hand, suppose  $\Gamma_i \vdash_{\omega K{\aclass{A}}} \varphi_i$. For every $j < i$, the above procedure cannot output $\langle \Gamma_i, \varphi_i\rangle$ since this has not been stored in $P$ yet. On the other hand, for every $j \geq i$, we know that the procedure to decide whether $\Gamma_i \vdash_{jK{\aclass{A}}} \varphi_i\rangle$ answers positively, and then the pair is never returned.
\end{proof}

The previous lemma implies that, for every class of algebras $\aclass{C}$ satisfying the premises of Theorem \ref{th:undec} for which
	$\models_{\aclass{C}}$ is decidable, the logic $\vdash_{\omega K{\aclass{C}}}$ is not recursively enumerable. Otherwise, since the previous lemma proves that $\{\langle \Gamma, \varphi \rangle \in \mathcal{P}_\omega(Fm) \times Fm \colon \Gamma \not \vdash_{\omega K{\aclass{C}}} \varphi\}$ is recursively enumerable, the logic $\vdash_{\omega K{\aclass{C}}}$ would be decidable, contradicting Theorem \ref{th:undec}. Since $\L$ and $\Pi$ are decidable logics \cite{Ha98}, their standard completeness implies that  $\models_{[0,1]_\sL}$ and $\models_{[0,1]_\Pi}$ are decidable too, leading to the following corollary.

	\begin{corollary}\label{cor:FinNotAxiomat}
	The logics $\vdash_{\omega K\sL}$ and $\vdash_{\omega K{\Pi}}$ are not axiomatizable.
\end{corollary}
	
	However, since it is not a general fact that the logics $\vdash_{K{\aclass{C}}}$ are complete with respect to finite models, the lack of axiomatization of the previous logics does still not close the problems mentioned in the beginning of this section.

\subsection{Modal \L ukasiewicz Logic is not axiomatizable}
We can show that, even if the global modal \L ukasiewicz logic might not enjoy the finite model property, if the global modal (standard) \L ukasiewicz logic $\vdash_{K\sL}$ is R.E. then $\vdash_{\omega K\sL}$ (the analogous logic over finite models) would be R.E. too. Since we saw before the latter does not hold, we will conclude that $\vdash_{K\sL}$ is not R.E. and so, not axiomatizable.

%

Lemma 3 from \cite{Ha05b} allows us to prove completeness of $\vdash_{\omega K\sL}$ with respect to witnessed models, in a similar way to how it is done for tautologies of fuzzy description logic (FDL)over \L ukasiewicz logic in the same publication. We do not introduce details of F.O. (standard) \L ukasiewicz logic here, we refer the interested reader to eg. \cite{Ha98}. Just recall that:  
\begin{itemize}
	\item A standard \L ukasiewicz F.O. model is a structure $\langle W, \{P_i\}_{i \in I}\rangle$ where $W$ is a non-empty set and for each $i \in I$ and $ar(i)$ the arity of $P_i$, $P_i \colon W^{ar(i)} \rightarrow [0,1]$, 
	\item An evaluation in a (F.O.) model is a mapping $v \colon \mathcal{V} \mapsto W$. Moreover, we write $v[x \mapsto m]$ to denote the evaluation $v$ where the mapping of the variable $x$ is overwritten and $x$ is mapped to $m$ (and simply $[x \mapsto m]$ denotes that the evaluation of $x$ is $m$ and the other variables are irrelevant).
	\item The value of a formula $\varphi$ in a (F.O.) model $\mod{M}$ under an evaluation $v$, denoted by $\|\varphi\|_{\mod{M}, v}$  is inductively defined by
	\begin{itemize}
		\item $\|P_i(x_1, \ldots, x_{ar(i)})\|_{\mod{M}, v} = P_i(v(x_1), \ldots, v(x_{ar(i)}))$;
		\item $\|\varphi_1 \star \varphi_2\|_{\mod{M}, v}  = \|\varphi_1\|_{\mod{M}, v} \star  \|\varphi_2\|_{\mod{M}, v}$ for $\star$ propositional (\L) operation;
		\item $\|\forall x \varphi(x)\|_{\mod{M}, v} = \bigwedge_{m \in W} \|\varphi\|_{\mod{M}, v[x \mapsto m]}$, 
		\item $\|\exists x \varphi(x)\|_{\mod{M}, v} = \bigvee_{m \in W} \|\varphi\|_{\mod{M}, v[x \mapsto m]}$.
	\end{itemize}
\end{itemize}
Observe that the value of a  sentence (i.e., closed formula, without free variables) $\varphi$ in a model is constant under any evaluation, so we can simply write $\|\varphi\|_{\mod{M}}$ to denote its value in a model. Moreover, we say that a model $\mod{M}$ is \textbf{witnessed} whenever for every sentence $Qx \varphi(x)$ for $Q \in \{\forall, \exists\}$ there is some $m \in W$ such that \[ \|Qx \varphi(x)\|_{\mod{M}} = \|\varphi(x)\|_{\mod{M}, [x \mapsto m]}.\]

The consequence relation over standard \L ukasiewicz F.O. models,  $\models_{\forall[0,1]_{\tiny{\L}}}$, is defined for sentences by letting $\Gamma  \models_{\forall[0,1]_{\tiny{\L}}} \varphi$ whenever for every standard \L ukasiewicz F.O. model $\mod{M}$, if $\|\Gamma\|_{\mod{M}} \subseteq \{1\}$ then $\|\varphi\|_{\mod{M}} = 1$.

\begin{lemma}[\cite{Ha05b}, Lemma 3]
	Let $\mod{M}$ be a standard \L ukasiewicz F.O. model. Then there is a (standard \L ukasiewicz F.O.) witnessed model $\mod{M}'$ such that $\mod{M}$ is a submodel of $\mod{M}'$ and for every sentence 
 	 $\alpha$ it holds that \[\|\varphi\|_{\mod{M}} = 1 \text{ if and only if } \|\varphi\|_{\mod{M}'} = 1.\]
\end{lemma}

From here, we can easily prove completeness of $\vdash_{K\sL}$ with respect to witnessed Kripke models, i.e., those for which, for every modal formula $\triangledown\varphi$ (with $\triangledown \in \{\Box, \Diamond\}$) and every world $v$ there is some world $w$ such that $Rvw$ and 
\[e(v, \triangledown\varphi) = (w, \varphi). \]
\begin{lemma}\label{lemma:witmodel}
	If $\Gamma \not \vdash_{K\sL} \varphi$ there is a witnessed standard \L ukasiewicz Kripke model $\mod{M}$ and $v \in W$ such that $\Gamma \not \vdash_{\langle \mod{M},v\rangle} \varphi$.
\end{lemma}
\begin{proof}
	We can use the usual translation from modal to F.O. logics in order to move from a Kripke model to a suitable F.O. model. 
	For a modal formula $\varphi$ consider the F.O. language $\{R/2, \{P_p/1\colon p\text{ variable in }\varphi\}\}$. For an arbitrary natural number $i \in \mathds{N}$, let us define the translation $\langle \varphi, x_i \rangle^*$ recursively by letting
	\begin{itemize}
		\item  $\langle p, x_i \rangle^* \coloneqq P_p(x_i)$;
		\item $\langle \varphi \star \psi, x_i \rangle^* \coloneqq \langle \varphi, x_i \rangle^* \star \langle \psi, x_i \rangle^*$ for a propositional connective $\star$;
		\item $\langle \Box \varphi, x_i \rangle^* \coloneqq \forall x_{i+1} R(x_i, x_{i+1}) \rightarrow \langle \varphi, x_{i+1} \rangle^*$;
		\item $\langle \Diamond \varphi, x_i \rangle^* \coloneqq \exists x_{i+1} R(x_i, x_{i+1}) \cdot \langle \varphi, x_{i+1} \rangle^*$;
	\end{itemize}

 It is a simple exercise to check that 
 \[\Gamma \vdash_{K\sL} \varphi \text{ if and only if } \{\forall x_0 \langle \gamma, x_0 \rangle^*\}_{\gamma \in \Gamma}, \forall x\forall y (R(x,y) \vee \neg R(x,y)) \models_{\forall[0,1]_{\tiny{\L}}} \forall x_0 \langle \varphi, x_0\rangle^*.\]

If $\Gamma \not \vdash_{K\sL} \varphi$ there is some F.O. model satisfying the premises of the right side of the above consequence and not $\forall x_0 \langle \varphi, x_0\rangle^*$. From the previous lemma we know there is a witnessed (F.O.) model $\mod{M}$ in which the same conditions hold. Then, there is some $m$ in the universe for which $\|\langle \varphi,x_0 \rangle^*\|_{\mod{M}, [x_0 \mapsto m]} < 1$. At this point, it is only necessary to build a witnessed Kripke model $\what{\mod{M}}$ from $\mod{M}$ that is a global model for $\Gamma$ but does not satisfy $\varphi$ at some world. In order to do that, let the universe of the Kripke model be the same universe of $\mod{M}$, and let the accessibility relation be given by the interpretation of the binary predicate $R$ in $\mod{M}$. Observe that, since $\forall x\forall y (R(x,y) \vee \neg R(x,y))$ is true in the model, necessarily $R(x,y) \in \{0,1\}$, and thus the resulting model will be crisp. Finally, let $e(v, p) = \|P_p(v)\|_{\mod{M}}$ for each variable $p$ and each world $v \in W$. 

By induction on the complexity of the formula it is routine to check that for every $\psi \in SFm(\Gamma, \varphi)$ and every $v \in W$, $e(v, \psi) = \|\langle \psi, x_i\rangle^*\|_{\mod{M}, [x_i \mapsto v]}$. Moreover, since the F.O. model is witnessed,  the Kripke model is witnessed too. $\what{\mod{M}}$ is a global model of $\Gamma$, while $e(m, \varphi) = \|\langle \varphi, x_0\rangle\|_{\mod{M}, [x_0 \mapsto m]} < 1$, concluding the proof of the lemma.
\end{proof}

We can use the non-idempotency of the \L ukasiewicz t-norm to recursively reduce the global consequence relation over finite models to the unrestricted global consequence relation.

\begin{lemma}\label{lemma:GLToFGL}
Let $\Gamma \cup \{\varphi\} \subseteq_{\omega} Fm$ and $p,q \not \in \mathcal{V}(\Gamma, \varphi)$. 
Define 
	\begin{itemize}
		\item $\Xi(p)\coloneqq \{\Box \const{0} \vee (p \leftrightarrow \Box p), \Box \const{0} \vee (\Box p \leftrightarrow \Diamond p)\}$,
		\item $\xi(p,q) \coloneqq (q \leftrightarrow p) \cdot \Box q$,
		\item $\psi(p,q) \coloneqq p \vee \neg p \vee q \vee \neg q$.
	\end{itemize}
Then

\[	\Gamma \vdash_{\omega K\sL} \varphi \text{ if and only if }\Gamma, \Xi(p), \xi(p,q) \vdash_{K\sL} \varphi \vee \psi(p,q)\] 
\end{lemma}
\begin{proof}
	$\Rightarrow$: Assume $\Gamma, \Xi(p), \xi(p,q) \not \vdash_{K\sL} \varphi \vee \psi(p,q)$. From Lemma \ref{lemma:witmodel} it follows that there is a witnessed standard \L ukasiewicz Kripke model $\mod{M}$ and $v \in W$ such that $\mod{M} \models \Gamma, \Xi(p), \xi(p,q)$ l,.
	(i.e. $e(u, \Gamma, \Xi(p) \cup \xi(p,q)) \subseteq \{1\}$ for all $u$) and $e(v, \varphi \vee \psi(p,q)) < 1$. We can assume that $\mod{M}$ is the unraveled tree generated from $v$. We will now prove that we can define a finite model equivalent to this one for what concerns the formulas in $F = \SF(\Gamma \cup \Xi(p) \cup \{\xi(p,q)\} \cup \{\varphi \vee \psi(p,q)\})$.
	
	First, since
	$e(u, \Xi(p)) \subseteq \{1\}$ for each $u \in W$, it follows that there is $a \in [0,1]$ such that for all $u \in W$, $e(u, p) = a$, as it was proven in Lemma \ref{lemma:valuey}. Moreover, from $e(v, p \vee \neg p) < 1$ we have that $a \in (0,1)$.

	On the other hand, from $e(u,  \xi(p,q)) =1$ for each $u \in W$ we get that  $e(u, q) = e(u, \Box q) a$. Thus, for each world $u \in W$ we have that $e(u, q) \leq a^s$ for all $s \leq \height(u)$.
	In particular, if there was any $u \in W$ with $\height(u) = \infty$, 
	$v$ would also have infinite height, and so $e(v, q) \leq a^n$ for all $n\in \mathds{N}$. Since $a \in (0,1)$, by the definition of the product in the standard MV algebra, the previous family of inequalities would imply that $e(v, q) = 0$. Then it would holds that $e(v, \neg q) = 1$, which is not possible since by the assumption
	$e(v, \varphi \vee \psi(p,q)) < 1$. Thus, $v$ -and so, all worlds of the model- must have finite height. 
	
	We can apply now a filtration-like transformation to $\mod{M}$ with respect to the set of formulas $F$ in order to obtain a finite directed model. 	
	To do this, let us denote by $wit(u, \triangledown \chi)$ an arbitrary witnessing world for the modal formula $\triangledown \chi$ at world $u$ (i.e., such that $e(u, \triangledown \chi) = e(wit(u, \triangledown \chi), \chi)$. Then define the universe $W' \coloneqq \bigcup_i \in \omega W_i$ with 
	\begin{align*}
	W_0 &\coloneqq \{v\}\\
	W_i+1 &\coloneqq \{wit(u, \triangledown \chi) \colon \triangledown \chi \in SFm(F), u \in W_i\}
	\end{align*}
	Observe that if there is some $i \in \mathds{N}$ for which the worlds in $W_i$ do not have successors, $W_j = \emptyset$ for any $j > i$. 
	
	We proved above that all worlds in $\mod{M}$ have finite height. Further, $F$ is a finite set of formulas. Henceforth, the model $\mod{M}'$ resulting from restring $\mod{M}$ to the universe $W'$ is a finite directed model with root $v$. Moreover, it is such that $e'(w, \Gamma, \Xi(p), \xi(p,q)) \subseteq \{1\}$ for each world $w \in W'$, and $e'(v, \varphi \vee \psi(p,q)) < 1$. In particular, $e(w, \Gamma) \subseteq \{1\}$ at each world $w$, and $e(v, \varphi) < 1$.
	
	$\Leftarrow$: Assume $\Gamma \not \vdash_{\omega K\sL} \varphi$, so there is a finite model $\mod{M}$ and a world $v \in W$ such that $\Gamma \not \vdash_{\langle \mod{M},v\rangle} \varphi$. 
	Let $\mod{M}'$ be the model with universe equal to that of $\mod{M}$, such that $e'(u,x) = e(u,x)$ for all $x \in \mathcal{V}(\Gamma, \varphi)$ and each $u \in W$, and where the values of variables $p,q$ are defined as follows. Pick  an arbitrary element $a \in (\frac{\height{v}}{\height{v}+1}, 1)$ and let, for each $u \in W$, 
	\begin{itemize}
		\item $e(u, p) = a$, 
		\item $e(u, q) = e(u, \Box q) a$ (observe this is well defined since all worlds have finite height, so we can define $q$ inductively from the worlds with height $0$).
	\end{itemize}
	This evaluation satisfies, in all worlds of the model, all formulas from $\Xi(p)$ and $\xi(p,q)$, and it forces $e(v, p) \not \in \{0,1\}$ and $e(v, q) \not \in \{0,1\}$. Moreover, it satisfies the formulas from $\Gamma$, and $e(v, \varphi) <1$, since the evaluation of all the variables appearing in $\Gamma$ and $\varphi$ has been preserved. Thus, $\Gamma, \Xi(p), \xi(p,q) \not \vdash_{K\sL} \varphi \vee \psi(p,q)$ either.
\end{proof}

The fact that the (finitary) \L ukasiewicz global modal logic is not axiomatizable follows as a consequence of previous reduction (which is recursive) and the undecidability of $\vdash_{\omega K\sL}$.

\begin{theorem}\label{th:LnotAx}
$\vdash_{K\sL}$ is not axiomatizable.
\end{theorem}
\begin{proof}
	Assume $\vdash_{K\sL}$ is axiomatizable, and so, recursively enumerable. 
	We can prove that then $\vdash_{\omega K\sL}$ is recursively enumerable too, contradicting Corollary \ref{cor:FinNotAxiomat}. For that, take a recursive enumeration of all pairs $\langle \Gamma, \varphi \rangle \in \mathcal{P}_\omega(Fm) \times Fm$ such that $\Gamma \vdash_{K\sL} \varphi$. 
	For each pair, let $\mathcal{V} = \mathcal{V}ars(\Gamma, \varphi)$, and check whether there are some $p,q \in \mathcal{V}$ for which $\Gamma = \Gamma_0(\mathcal{V}\setminus\{p,q\}) \cup \Xi(p) \cup \{\xi(p,q)\}$ and 	$\varphi = \varphi_0(\mathcal{V}\setminus\{p,q\}) \vee \psi(p,q)$ for some $\Gamma_0, \varphi_0$. This is a decidable procedure because $\Gamma$ is a finite set and the translations $\Gamma_0$ and $\varphi_0$ are recursive. 
If that is the case, output $\langle \Gamma_0, \varphi_0\rangle$, and don't output anything otherwise.
Lemma \ref{lemma:GLToFGL} implies that this procedure enumerates $\vdash_{\omega K\sL}$.

	However, Corollary \ref{cor:FinNotAxiomat} states that $\vdash_{\omega K\sL}$ is not R.E., a contradiction.
\end{proof}

\subsection{Modal Product Logic is not axiomatizable either} 
In \cite{BaHaKraSve98}, the authors observe, roughly speaking, that the standard MV algebra is isomorphic to the standard product algebra restricted to $[a,1]$ for arbitrary fixed $0 < a < 1$.
Relying on the isomorphism there provided, 
they also show that the tautologies of standard \L ukasiewicz propositional and F.O. logics\footnote{Standard \L ukasiewicz propositional logic is indeed the \L ukasiewicz propositional logic. However, the F.O. logic over the standard MV algebra and the analogous logic over all chains in the variety differ \cite{Ha98}.} can be recursively reduced to those of the respective (standard) Product logic. In \cite[Lem. 4.1.14, Lem. 6.3.5]{Ha98} these results are formulated regarding the corresponding logical deduction relations.

We can use a similar argument, slightly modifying the reduction so it works in the modal case.\footnote{It is possible to do an alternative proof reducing $\vdash_{\omega K\Pi}$ to $\vdash_{K\Pi}$, similar to the one in the previous section.}
 
Given a finite set of variables $\mathcal{V}$, let $x$ be a propositional variable not in $\mathcal{V}$. For each formula 
$\varphi$ of $K\sL$ in variables $\mathcal{V}$, define its translation $\varphi^x$ as follows:
\begin{align*}
(0)^x \coloneqq& x & (q)^x \coloneqq& q \vee x \text{ for each }q \neq x,\\
(\varphi \rightarrow \psi)^x \coloneqq& (\varphi^x \rightarrow \psi^x) & (\varphi \cdot \psi)^x\coloneqq& x \vee (\varphi^x \cdot \psi^x)\\
(\Box \varphi)^x\coloneqq & \Box \varphi^x& &
\end{align*}

Further, let $\Theta^x \coloneqq  \{\Box x \leftrightarrow \Diamond x, \Box x \leftrightarrow x, \neg \neg x\}$. 

In the spirit of Lemmas 2 and 3 from \cite{BaHaKraSve98} it is possible to prove the following result. The proof is very similar to the one in the previous reference, but for the sake of completeness we provide the details in the appendix.
\begin{lemma}\label{lem:redluktoprod}
For every $x \not \in \mathcal{V}ar(\Gamma \cup \{\varphi\})$, it holds that $\Gamma \vdash_{K\sL} \varphi$ if and only if $\Gamma^x, \Theta^x \vdash_{K\Pi} \varphi^x$.
\end{lemma}

Since the reduction is recursive, together with Theorem \ref{th:LnotAx}, the following is immediate.
\begin{corollary}\label{Cor:PnotAx}
$\vdash_{K\Pi}$ is not axiomatizable.
\end{corollary}

We have proven that $\vdash_{K\sL}$ and $\vdash_{K\Pi}$ are not in $\Sigma_1$ from the Arithmetical Hierarchy. 
We leave open the question of whether they are $\Pi_2$-complete, as it is the case for the tautologies of their F.O. versions, or whether they belong to some other level of the hierarchy. The proofs of Ragaz in \cite{Ra83, Ra83c}  heavily rely on the  expressive power of F.O. logic, and also in proving the result directly for tautologies of the logic. In the present work, the results affect the logic itself, since for instance, the tautologies of $\vdash_{K\sL}$ are decidable: they coincide with those of $\vdash^l_{K\sL}$ and this logic is decidable (\cite[Corollary 4.5]{Vi20}).

\section{The necessitation rule}\label{sec:necRule} 

Recall that in (classical) modal logic, the global deduction is axiomatized as the local one plus the (unrestricted) necessitation rule $N_\Box\colon x \vdash \Box x$. It was asked in \cite{BoEsGoRo11} whether this was the case in general, or whether at least, this condition held for modal expansions of fuzzy logics. This was the case in the modal logics known up to now (eg. in modal G\"odel logics, and in the infinitary modal \L ukasiewicz and Product logics studied in the literature). The question remained open in full generality over modal many-valued logics. We can give a negative answer to the problem, first by a simple counter-example over the modal expansions of \L ukasiewicz logic (using the non-axiomatizability of $\vdash_{K_\sL}$ proven in Theorem \ref{th:LnotAx}) and later proving that this is a more general fact.

First, it is possible to see that the local deduction is decidable using the version of Lemma \ref{lemma:witmodel} referring to the local logic. A detailed proof of the decidability of $\vdash^l_{K_\sL}$ can be found in \cite[Corollary 4.5]{Vi20}. Thus, $\vdash^l_{K_\sL}$ has a recursive axiomatization (for instance, built by enumerating all possible pairs $\langle \Gamma, \varphi \rangle$ with $\Gamma \cup \{ \varphi\} \subseteq_{\omega} Fm$, and then returning the pairs for which $\Gamma \vdash^l_{K_\sL} \varphi$). On the other hand, if the global consequence were to coincide with the local one plus the $N_\Box$ rule, the logic axiomatized by adding to the previous system the $N_\Box$ rule should produce a recursive axiomatization of $\vdash_{K_\sL}$, contradicting Theorem  \ref{th:LnotAx}. 

As we said, it is possible to widen the scope of the previous result, and produce a constructive proof serving all modal logics built over classes of algebras like the ones in Theorem \ref{th:undec}. This can be done following an approach different from the previous direct one working for the \L ukasiewicz case, and instead providing a derivation that is valid in the global modal logics and not in the local ones extended by the necessitation rule.

 For simplicity, let us fix a class of algebras $\aclass{A}$ like the one from Section \ref{sec:undecGlobal}, and let $\vdash$ and $\vdash^l$ denote $\vdash_{K_{\aclass{A}}}$ and $\vdash^l_{K_{\aclass{A}}}$ respectively. Further, let $\vdash^l_{N_\Box}$ denote the logic $\vdash^l$ plus the necessitation rule $x \vdash \Box x$. 
A natural way to understand this extension is by considering the (possibly non recursive) list of finite derivations valid in $\vdash^l_{K_{\aclass{A}}}$ and add to this set the rule schemata $N_\Box$. Let us call this set $R$. The minimal logic  containing $R$, namely $R^l$, is the logic $\vdash^l_{N_\Box}$. Since all rules in $R$ have finitely many premises, the resulting logic is finitary.

Considering $R$ as a (possibly non recursive) axiomatization for $\vdash^l_{N_\Box}$, all derivations valid in $\vdash^l$ have a proof in the extended system of length $0$. Thus, the length of the proofs in the extended system only reflects the applications of the necessitation rule. Since by definition $\vdash^l$ is a finitary logic, only finitely many applications of the rule are used at each specific derivation.
This means that a proof of $\varphi$ from a finite set of premises $\Gamma$ in this axiomatic system is given simply as a finite list of pairs $\langle \Gamma_i, \varphi_i\rangle_{0 \leq i \leq N}$ such that 
\begin{itemize}
\item $\Gamma_0 = \Gamma$ and $\varphi_N = \varphi$, 
\item For each $0 \leq i \leq N$, $\Gamma_i \vdash^l \varphi_i$ and, 
\item $\Gamma_{i+1} = \Gamma_i \cup \{\Box \varphi_i\}$. 
\end{itemize}

From  here, it is quite simple to prove the following characterization of $\vdash^l_{N_\Box}$.
\begin{lemma}
	\[\Gamma \vdash^l_{N_\Box} \varphi \text{ if and only if } \{\Box^i\Gamma\}_{i \in \mathds{N}}\vdash^l \varphi.\]
\end{lemma}
\begin{proof}
The right to left direction is immediate. For the other direction, if $\Gamma \vdash^l_{N_\Box} \varphi$, since the logic is finitary, $\varphi$ can be proven from $\Gamma$ by using the $N_\Box$ rule a finite number of times, say $n$.
It can be easily proven by induction in $n$ \footnote{Using that $\Gamma \vdash^l \psi \Rightarrow \Box \Gamma \vdash^l \Box \psi$.} that $\{\Box^i\Sigma \}_{i \leq n} \vdash^l \chi$ if and only if $\Sigma \vdash^l_{n \cdot N_\Box} \chi$, where $n\cdot N_\Box$ stands for using  the $N_\Box$ rule up to $n$ times.  
This concludes the proof.
\end{proof}

We can then produce a set of formulas that yields a valid derivation in the global logics, but it does not in the corresponding local logics plus necessitation.
\begin{theorem}
	$\vdash$ does not coincide with $\vdash^l_{N_\Box}$.
\end{theorem}
\begin{proof}
We claim that both
\begin{align*}
y  \leftrightarrow \Box y, y \leftrightarrow \Diamond y, x \leftrightarrow (\Box x) y, \neg \Box \perp & \vdash x \rightarrow xy, \text{ and }\\
y  \leftrightarrow \Box y, y \leftrightarrow \Diamond y, x \leftrightarrow (\Box x) y, \neg \Box \perp & \not \vdash^l_{N_\Box} x \rightarrow xy
\end{align*}
which proves the theorem.

Regarding the first claim, consider an arbitrary Kripke model satisfying globally the set of premises. In particular, from $\neg \Box \perp$ we get that each world in the model has a successor, and so, has infinite height in the sense of Definition \ref{def:height}. Moreover, the value of $y$ is constant inside each connected part of the model, as in Lemma \ref{lemma:valuey}. Consider each connected submodel $\mod{M}$, and let $\alpha$ be the value of $y$ in it. Then, at each point of the model, $e(u, x) \leq \alpha^i$ for all $i\in \mathds{N}$. Then, since the algebras in the class are weakly saturated, we get that $e(u,x) \alpha = e(u,x)$, proving the formula in the right side.

In order to prove the second claim, let us denote by  $\Sigma$ the set of premises. From the previous lemma we have that our claim holds if and only if 
\[\{\Box^i \Sigma\}_{i \in \mathds{N}} \not \vdash^l  x \rightarrow xy.\] 
Since $\vdash^l$ is finitary by definition, this holds if and only if $\{\Box^i\Gamma\}_{i \leq N} \not \vdash^l \varphi$ for all  $N \in \mathds{N}$. We can produce a counter-model for each $N \in \mathds{N}$.

Indeed, consider a model with universe $\{0, \ldots, N+1\}$, and the accessibility relation given by $R = \{\langle i, i+1\rangle \colon i \leq N\}$. Regarding the evaluation, pick an arbitrary $\alg{A} \in \aclass{A}$ that is not $(N+1)$-contractive, and pick $a \in \alg{A}$ such that $a^{N+2} < a^{N+1}$. Then let
\[e(i, y) = a \text{ for } 1 \leq i \leq N+1, \qquad e(N+1, x) = 1, \qquad e(i, x) = a^{N+1-i} \text{ for } 1 \leq i \leq N.\]

It is routine to check that this evaluation satisfies $\{\Box^i\Gamma\}_{i \leq N} $ at the world $0$, i.e., $e(0, \Box^i\Sigma) = 1$ for all $i \leq N$. On the other hand, observe that $e(0, x) = a^{N+1}$. Due to the way  we chose $a$ it holds that $e(0, x y)  = a^{N+2} < a^{N+1} = e(0,x)$, thus falsifying the consequence. 

\end{proof}

\noindent
\textbf{Acknowledgments.} The author is thankful to the anonymous referee for many useful and detailed contributions. This project has received funding from the following sources: the European Union’s Horizon 2020 research and innovation programme under the Marie Skłodowska-Curie grant agreement No. 101027914; the Grant No. CZ.02.2.69/0.0/0.0/17 050/0008361 of the Operational programme Research, Development, Education of the Ministry of Education, Youth and Sport of the Czech Republic. I also wish to thank Lluís Godo, Gavin St. John and Tommaso Moraschini for their useful help.

\bibliographystyle{plainurl}

\newpage
\section{Appendix}
\subsection*{Craig's Theorem}
It is easy to see that a finitary R.E. logic with a definable idempotent $n$-ary operation for every $n$ is always axiomatizable. This can be checked as it is done in Craig's Theorem for classical logic.

Indeed, let $\star_n$ be an idempotent n-ary operation (i.e., such that $\star_n(x,\ldots,x) \vdash_{\mathcal{L}} x$ and $x \vdash_{\mathcal{L}} \star_n(x,\ldots, x)$), as for instance, the classical conjunction $\wedge$ in classical logic. For convenience, denote $\star_n(\varphi, \ldots, \varphi)$ by $\star_n(\varphi)$. 
Consider a recursive enumeration of the finitary derivations of the logic, i.e.,  $\{\langle \Gamma_n, \varphi_n\rangle \colon n \in \omega\} = \mathcal{L}^{fin}$ (which exists because the logic is R.E), and take the axiomatization given by the set of rules
\[R = \{\langle \Gamma_n, \star_n(\varphi_n)\rangle, \langle \star_n(\varphi_n), \varphi_n \rangle\colon n \in \omega\}\]
This is a recursive set: for a finitary rule $\langle \Sigma, \chi\rangle$, it is first possible to decide whether there is some $1 \leq j \leq symb(\chi)$ (where $symb(\chi)$ denotes the number of symbols appearing in $\chi$) such that 
$\chi$ is of the form $\star_j(\varphi_j)$ and $\Sigma = \Gamma_j$.
If this is the case,  $\langle \Sigma, \chi\rangle$ belongs to $R$. Otherwise, it is possible to decide whether $\Sigma$ is a singleton $\{\psi\}$ and, in this case, whether there is some $1 \leq j \leq symb(\psi)$ such that $\psi$ is of the form $\star_j(\varphi_j)$ and  $\chi$ is exactly $\varphi_j$. If the answer to both questions is yes, $\langle \Sigma, \chi\rangle$ belongs to $R$. Otherwise, the rule does not belong to $R$.

Further, it is clear that $\mathcal{L}\subseteq R^l$. On the other hand, since the rules in $R$ are finitary, it is easy to see (by means of the characterization of $R^l$ as the proofs in $R$) that $R^l \subseteq \mathcal{L}$ too.

\subsection*{Proof of Lemma \ref{lem:redluktoprod}}

The proof of Lemma \ref{lem:redluktoprod} detailed below draws inspiration from the results in \cite{BaHaKraSve98}, and relies in the same isomorphic mappings introduced there. However, the approach and details are slightly different here, since we formulate alternative intermediate results, and we propose a more explicit proof using basic arithmetics.

\textit{Proof of Lemma \ref{lem:redluktoprod}: (For every $x \not \in \mathcal{V}ar(\Gamma \cup \{\varphi\})$, it holds that $\Gamma \vdash_{K\sL} \varphi$ if and only if $\Gamma^x, \Theta^x \vdash_{K\Pi} \varphi^x$.)}\\
$\Leftarrow$: Assume $\Gamma \not \vdash_{K\sL} \varphi$. Then there is some standard \L ukasiewicz Kripke model $\mod{M}$ such that $\mod{M} \models \Gamma$ but $\mod{M},v \not \models \varphi$ for some $v$ in the model. Chose an arbitrary $a \in (0,1)$, and let us define a standard product model $\mod{M}'$ by letting the universe and accessibility relations be those of $\mod{M}$, and further, for each $w \in W$, let\footnote{The mapping below is the isomorphism between the standard MV algebra and the product algebra restricted to $[a,1]$ used in \cite{BaHaKraSve98}.}
\begin{align*}
e'(w, x) \coloneqq& a & e'(w, q) \coloneqq& a^{1-e(w,q)} \text{ for each variable }q \neq x.
\end{align*}

\textbf{Claim 1:}\textit{For every formula $\psi$ in variables from $\mathcal{V}ar(\Gamma \cup \{\varphi\})$, and for each $w \in W$, it holds that 
\[e'(w, \psi^x) = a^{1-e(w, \psi)}\]}
The conclusion easily follows, as we proceed to explain. Indeed, Claim 1 implies that $e'(w, \gamma^x) = a^{1-e(w, \gamma)} = a^0 = 1$ for each $\gamma \in \Gamma$ and $w \in W$. It is also clear that $e'(w, \Theta^x) = 1$, since $x$ is evaluated to the same element $a > 0$ in all worlds of $\mod{M}'$. On the other hand, $e'(v, \varphi^x) = a^{1-e(v, \varphi)} <1$, since $e(v, \varphi) <1$. Thereby, $\mod{M}' \models \Gamma^x, \Theta^x$ and $\mod{M}',v \not \models \varphi^x$, and so, $\Gamma^x, \Theta^x \not \vdash_{K\Pi} \varphi ^x$. Then, it only remains to prove the claim.

\textbf{Proof of Claim 1. }
It can be proven by by induction on the complexity of the formula.
\begin{itemize}
\item For variables it is straightforward from the definition, since $e'(w,q^x) = e'(w,q \vee x) = a \vee a^{1-e(w,q)}$, and $a \leq a^q$ for every $q \in [0,1]$. 
\item For $\psi = \psi_1 \cdot \psi_2$, we have the following chain of equalities 
\begin{align*}
&\ e'(w, (\psi_1 \cdot \psi_2)^x) =\ e'(w, x \vee (\psi_1^x \cdot \psi_2^x)) =\ e'(w, x) \vee (e(w, \psi_1^x) \cdot_\Pi e(w, \psi_2^x)) \\
 \overset{I.H}{=}&\ a \vee (a^{1-e(w, \psi_1)} \cdot_\Pi a^{1-e(w, \psi_2)}) =\ a^{1 - (e(w, \psi_1) + e(w, \psi_2) -1)} 
= a^{1 - (e(w, \psi_1) \cdot_\sL e(w, \psi_2))} \\
 =&\ a^{1 - e(w, \psi_1 \cdot \psi_2)} = a^{1 -e(w, \psi)}.
\end{align*}

\item For $\psi = \psi_1 \rightarrow \psi_2$, we have the following chain of equalities 
\begin{align*}
& e'(w, (\psi_1 \rightarrow \psi_2)^x) = e'(w, \psi_1^x \rightarrow \psi_2^x) 
\overset{I.H}{=} a^{1 - e(w, \psi_1)} \rightarrow_\Pi a^{1 - e(w, \psi_2)} \\
= &\ \begin{cases} 1 &\hbox{ if } a^{1 - e(w, \psi_1)}  \leq a^{1 - e(w, \psi_2)} \\ \frac{a^{1 - e(w, \psi_2)}}{a^{1 - e(w, \psi_1)}} & \hbox{ otherwise} \end{cases} 
= \begin{cases} 1 &\hbox{ if } e(w, \psi_1) \leq e(w, \psi_2) \\ a^{1 - e(w, \psi_2) - 1 + e(w, \psi_1)} & \hbox{ otherwise} \end{cases}\\
=&\ \begin{cases} 1 &\hbox{ if } e(w, \psi_1)  \leq e(w, \psi_2) \\ a^{1 - (e(w, \psi_1) \rightarrow_\sL e(w, \psi_2)} & \hbox{ otherwise} \end{cases}\\
=&\  a^{1 - e(w, \psi_1 \rightarrow \psi_2)}.
\end{align*}

\item For $\psi = \Box \psi_1$, we know that
\begin{align*}
e'(w, (\Box \psi_1)^x) & = e'(w, \Box \psi_1^x) =  \bigwedge_{Rwu} e'(u, \psi_1^x) \overset{I.H}{=} \bigwedge_{Rwu} a^{1-e(u, \psi_1)}
\end{align*}
On the one hand, since $e(w, \Box \psi_1) = \bigwedge_{Rwu}e(u, \psi_1) \leq e(u, \psi_1)$ for each $u$ with $Rwu$, it holds that $a^{1-e(w, \Box \psi_1)} \leq \bigwedge_{Rwu} a^{1-e(u, \psi_1)} = e'(w, \Box \psi_1)$.

On the other hand,  
\begin{align*}
& \bigwedge_{Rwu} a^{1-e(u, \psi_1)} \leq a^{1-e(u, \psi_1)}\ \forall u \text{ s.t }Rwu
 & \Longrightarrow &\ a^{e(u, \psi_1)} \leq \frac{ a}{\bigwedge_{Rwu} a^{1-e(u, \psi_1)}}\ \forall u \text{ s.t } Rwu\\
 \overset{a \in (0,1)}{\Longrightarrow} & e(u, \psi_1) \geq log_a(\frac{ a}{\bigwedge_{Rwu}a^{1-e(u, \psi_1)}})\  \forall u \text{ s.t }Rwu
& \Longrightarrow & \bigwedge_{Rwu} e(u, \psi_1)  \geq log_a(\frac{ a}{\bigwedge_{Rwu}a^{1-e(u, \psi_1)}})\\
 \Longrightarrow &\ a^{\bigwedge_{Rwu} e(u, \psi_1)} \leq \frac{ a}{\bigwedge_{Rwu}a^{1-e(u, \psi_1)}}
 & \Longrightarrow & \bigwedge_{Rwu}a^{1-e(u, \psi_1)} \leq \frac{a}{a^{\bigwedge_{Rwu} e(u, \psi_1)}} \\
 \Longrightarrow &\ e'(w, (\Box \psi_1)^x) \leq a^{1 - \bigwedge_{Rwu} e(u, \psi_1)}
\end{align*}
This concludes the proof of the claim.
\end{itemize}

\smallskip

$\Rightarrow$: It is proven similarly, using the corresponding inverse of the isomorphism from \cite{BaHaKraSve98}. Assume there is a standard product model $\mod{P}$ such that $\mod{P} \models \Gamma^x, \Theta^x$ and $\mod{P},v \not \models \varphi^x$ for some $v$ in the model. As proven in Lemma \ref{lemma:valuey}, there is some element $a$ such that $e(w,x) = a$ for all $w$ in the universe of the model. Moreover, since $e(v, \neg \neg x) = 1$, necessarily $a > 0$.

Observe first that if $a = 1$ then $e(w,\psi^x) = 1$ for every $\psi$ with variables in $\mathcal{V}ar(\Gamma \cup \{\varphi\})$ and each $w$ in the universe. This observation is immediate by induction on the complexity of the formula. Since this would contradict the fact that $e(v, \varphi^x) < 1$, necessarily $a < 1$.

Let us define a standard \L ukasiewicz model $\mod{P}'$ as the model whose universe and accessibility relation are those of $\mod{P}$ and for each variable $q$ and each world $w$ let $e'(w, q) \coloneqq 1-log_a a \vee e(w, q)$. 

\textbf{Claim 2. }\textit{For every formula $\psi$ in variables from $\mathcal{V}ar(\Gamma \cup \{\varphi\})$ and for each $w \in W$ it holds that 
\[e'(w, \psi) = 1-log_a e(w, \psi^x)\]}

As in the previous case, the conclusion easily follows. Indeed, If Claim 2 holds, $e'(w, \gamma) = 1-log_a e(w, \gamma^x) = 1 - log_a 1 = 1$ for each $\gamma \in \Gamma$ and $w \in W$, and $e'(v, \varphi) = 1-log_a e(v, \varphi^x) = 1 - log_a \alpha$ for some $a \leq \alpha < 1$. Since the logarithm in base $a$ of elements in that interval is a value in $(0,1]$, necessarily $e'(v, \varphi)  < 1$, concluding the proof of the lemma. Then, it only remains to prove the claim.

\textbf{Proof of Claim 2.}

We will prove it by induction on the complexity of the formula. 
Observe a consequence of Claim 2 is that 
$e(w, \psi^x) \geq a$ for all $\psi$ and $w$ as before\footnote{Since $e'$ is defined inductively from the value of propostional variables, it always returns a value in $[0,1]$. Thus, for any $\psi$ the claim implies that $0 \leq 1-log_a e(w, \psi^x)$, so $log_a e(w, \psi^x) \leq 1$. This is only possible if $e(w, \psi^x) \geq a$.}. We will (inductively) use this property in the modal step, and refer to it by $I.H'$.

\begin{itemize}
\item For variables it is immediate, since $e'(w, q) =  1-log_a a \vee e(w, q) = 1-log_a e(w, q^x)$, 
\item For $\psi = \psi_1 \cdot \psi_2$, we have the following chain of equalities:\\
\noindent
$ e'(w, \psi_1 \cdot \psi_2) = e'(w, \psi_1) \cdot_\sL e'(w, \psi_2) \overset{I.H}{=} max\{0, 1-log_a e(w, \psi_1^x) + 1- log_a e(w, \psi_2^x) - 1\}
= max\{0,  1 - log_a(e(w,  \psi_1^x \cdot \psi_2^x))\}$.\\
\noindent
Now, if $0 < e(w, \psi_1^x \cdot \psi_2^x) < a$, it holds that $ log_a(e(w,  \psi_1^x \cdot \psi_2^x)) > 1$, and thus, $max\{0,  1 - log_a(e(w,  \psi_1^x \cdot \psi_2^x))\} = 0 = 1 - log_a(e(w,  \psi_1^x \cdot \psi_2^x)\vee a)$. It follows that 
\[max\{0,  1 - log_a(e(w,  \psi_1^x \cdot \psi_2^x))\} = 1-log_a(e(w,  \psi_1^x \cdot \psi_2^x)\vee a) = 1 - log_a e(w, (\psi_1 \cdot \psi_2)^x).\]

\item For $\psi = \psi_1 \rightarrow \psi_2$, we have the following chain of equalities:
\begin{align*}
&\ e'(w, \psi_1 \rightarrow \psi_2) = min\{1, 1 – e'(w, \psi_1) + e'(w, \psi_2)\}\\
 \overset{I.H}{=}&\ min\{1, 1 - (1-log_a e(w, \psi_1^x)) + 1- log_a e(w, \psi_2^x)\}\\
=&\  min\{1, 1 - (log_a e(w, \psi_2^x) - log_a e(w, \psi_1^x))\}\\
=&\  min\{1, 1 - log_a e(w, \psi_2^x)/e(w, \psi_1^x)\} \\
=&\ 1-log_a(min\{1, e(w, \psi_2^x)/e(w, \psi_1^x)\})\\
=&\ 1 - log_a e(w, \psi_1^x \rightarrow \psi_2^x)\\
 = &\ 1 - log_a e(w, (\psi_1 \rightarrow \psi_2)^x).
\end{align*}

\item For $\psi = \Box \psi_1$, we have that $e'(w, \Box \psi_1) = \bigwedge_{Rwv}e'(v,\psi_1) \overset{I.H}{=} \bigwedge_{Rwv} (1-log_a e(v, \psi_1^x)) = 1 - \bigvee_{Rwv}log_a e(v, \psi_1^x)$. Now, by $I.H'$, we know that $e(v, \psi_1^x) \geq a$ for all $v$ in the model, so in particular $\bigwedge_{Rwv} e(v, \psi_1^x) \in [a,1]$. Since the function $log_a()$ is continuous and decreasing in $[a, 1]$, it follows that $1 - \bigvee_{Rwv}log_a e(v, \psi_1^x) = 1 - log_a \bigwedge_{Rwv}  e(v, \psi_1^x) = 1 - log_a e(w, \Box \psi_1^x) = 1 - log_a e(w, (\Box \psi_1)^x)$. \qed

\end{itemize}

\end{document}